\documentclass[11pt]{article}
\usepackage{amssymb}
\begin{document}
\setlength{\oddsidemargin}{0cm}
\setlength{\evensidemargin}{0cm}
\baselineskip=20pt

\begin{center} {\Large \bf
Left-symmetric Bialgebras and An Analogue of the Classical
Yang-Baxter Equation}  \end{center}

\bigskip

\begin{center}  { \large Chengming ${\rm Bai}^{1,2}$} \end{center}

\begin{center}{\it 1. Chern Institute of Mathematics \& LPMC, Nankai University,  Tianjin 300071, P.R. China} \end{center}

\begin{center}{\it 2. Department of Mathematics, Rutgers, The State University of
New Jersey, Piscataway, NJ 08854, U.S.A.}\end{center}

\vspace{0.3cm}

\begin{center} {\large\bf   Abstract } \end{center}

We introduce a notion of left-symmetric bialgebra which is an
analogue of the notion of Lie bialgebra. We prove that a
left-symmetric bialgebra is equivalent to a symplectic Lie algebra
with a decomposition into a direct sum of the underlying vector
spaces of two Lagrangian subalgebras. The latter is called a
parak\"ahler Lie algebra or a phase space of a Lie algebra in
mathematical physics. We introduce and study coboundary
left-symmetric bialgebras and our study leads to what we call
``$S$-equation", which is an analogue of the classical Yang-Baxter
equation. In a certain sense, the $S$-equation associated to a
left-symmetric algebra reveals the left-symmetry of the products. We
show that a symmetric solution of the $S$-equation gives a
parak\"ahler Lie algebra. We also show that such a solution
corresponds to the symmetric part of a certain operator called
``${\cal O}$-operator", whereas a skew-symmetric solution of the
classical Yang-Baxter equation corresponds to the skew-symmetric
part of an ${\cal O}$-operator. Thus a method to construct symmetric
solutions of the $S$-equation (hence parak\"ahler Lie algebras) from
${\cal O}$-operators is provided. Moreover, by comparing
left-symmetric bialgebras and Lie bialgebras, we observe that there
is a clear analogue between them and, in particular, parak\"ahler
Lie groups correspond to Poisson-Lie groups in this sense.

\vspace {0.2cm}


\vspace{0.2cm}


\newpage

\section{Introduction}

Left-symmetric algebras (or under other names like pre-Lie algebras,
quasi-associative algebras, Vinberg algebras and so on) are
Lie-admissible algebras ((nonassociative) algebras whose commutators
are Lie algebras) whose left multiplication operators form a Lie
algebra. They have already been introduced by A. Cayley in 1896 as a
kind of rooted tree algebras ([Ca]). They also arose from the study
of several topics in geometry and algebra in 1960s, such as convex
homogenous cones ([V]), affine manifolds and affine structures on
Lie groups ([Ko],[Mat]), deformation of associative algebras ([G])
and so on. In particular, a Lie algebra ${\cal G}$ with a compatible
left-symmetric algebra structure is the Lie algebra of a Lie group
$G$ with a left-invariant affine structure, that is, there exists a
left-invariant (locally) flat and torsion free connection $\nabla$
in $G$. The left-symmetric algebra structure corresponds to the
connection $\nabla$ given by $XY=\nabla_XY$ for $X,Y\in {\cal G}$
and (a geometric interpretation of) the left-symmetry is just the
flatness of the connection $\nabla$ ([Me],[Ki]).

Furthermore, as it was pointed out in [CL] by Chapoton and Livernet,
the left-symmetric algebra ``deserves more attention than it has
been given''. It appears in many fields in mathematics and
mathematical physics. In [Bu2], Burde gave a survey of certain
different fields in which left-symmetric algebras play an important
role, such as vector fields, rooted tree algebras, words in two
letters, vertex algebras, operad theory, deformation complexes of
algebras, convex homogeneous cones, affine manifolds, left-invariant
affine structures on Lie groups (see [Bu2] and the references
therein). Here are some more examples (partly overlap with some
examples in [Bu2]).

(a) Symplectic structures on Lie groups and Lie algebras. A
symplectic Lie group is a Lie group $G$ with a left-invariant
symplectic form $\omega^+$. One can define an affine structure on
$G$ by ([Ch])
$$\omega^+(\nabla_{x^+}y^+,z^+)=-\omega^+(y^+,[x^+,z^+])\eqno (1.1)$$
for any left-invariant vector fields $x^+,y^+,z^+$ and hence
$x^+y^+=\nabla_{x^+}y^+$ gives a left-symmetric algebra. In fact,
equation (1.1) is of great importance to the study of symplectic and
K\"ahler Lie groups ([H],[LM],[DaM1-2],[MS],[KGM]).

(b) Complex and complex product structures on Lie groups and Lie
algebras. From a real left-symmetric algebra $A$, it is natural to
define a Lie algebra structure on the vector space $A\oplus A$ (that
is ${\cal G}(A)\ltimes_L {\cal G}(A)$) such that
$$J(x,y)=(-y,x),\;\;\forall\; x,y\in A\eqno (1.2)$$
is a complex structure on it. Moreover, there is a correspondence
between left-symmetric algebras and complex product structures on
Lie algebras ([AS]), which plays an important role in the theory of
hypercomplex and hypersympletic manifolds ([Bar],[AD]).

(c) Vertex algebras. Vertex algebras are fundamental algebraic
structures in conformal field theory ([FLM],[FHL],[JL]). For any
vertex algebra $V$,
$$a*b=a_{-1}b,\;\;\forall\; a, b\in V\eqno (1.3)$$
defines a left-symmetric algebra. And a vertex algebra is equivalent
to a left-symmetric algebra and a Lie conformal algebra with some
compatibility conditions ([BK]).

Vertex algebras are also closely related  to a subclass of
(finite-dimensional) left-symmetric algebras, namely, Novikov
algebras.  Novikov algebras are left-symmetric algebras with
commutative right multiplication operators. They were introduced in
connection with Hamiltonian operators in the formal variational
calculus ([GD]) and the following Poisson brackets of hydrodynamic
type ([BN])
$$\{ u(x_1), v(x_2)\}=\frac{\partial}{\partial x_1}
((uv)(x_1))x_1^{-1}\delta
(\frac{x_1}{x_2})+(uv+vu)(x_1)\frac{\partial}{\partial x_1}
x_1^{-1}\delta (\frac{x_1}{x_2}).\eqno (1.4)$$ Furthermore, let $A$
be a Novikov algebra and set ${\cal A}=A\otimes {\bf C}[t,t^{-1}]$,
where $t$ is an indeterminate.  Then the bracket
$$[a\otimes t^m, b\otimes t^n]
=(-mab+nba)\otimes t^{m+n-1},\;\;\forall a,b\in A,\;m,n\in {\bf Z}
\eqno (1.5)$$ defines a Lie algebra structure on ${\cal A}$ and this
Lie algebra can be used to construct a vertex Lie algebra and a
vertex algebra ([Li]). Conversely, vertex algebras satisfying
certain conditions must correspond to some Novikov algebras (roughly
speaking, such a vertex algebra $V$ is generated from $V_{(2)}$
which is a Novikov algebra, with some additional conditions)
([BKL]).

(d) Phase spaces of Lie algebras. The concept of phase space of a
Lie algebra was introduced by Kupershmidt in [Ku1] by replacing the
underlying vector space with a Lie algebra and was generalized in
[Bai2]. In [Ku2], Kupershmidt pointed out that left-symmetric
algebras appear as an underlying structure of those Lie algebras
that possess a phase space and thus they form a natural category
from the point of view of classical and quantum mechanics.

(e) Left-symmetric algebras are closely related to certain
integrable systems ([Bo1], [SS],[W]), classical and quantum
Yang-Baxter equation ([ES],[Ku3],[GS],[DiM]), combinatorics ([E])
and so on. In particular, they play a crucial role in the Hopf
algebraic approach of Connes and Kreimer to renormalization of
perturbative quantum field theory ([CK]).

In this paper, we study a structure, namely, parak\"ahler structure,
which appears in both geometry and mathematical physics, in terms of
left-symmetric algebras. In geometry, a parak\"ahler manifold is a
symplectic manifold with a pair of transversal Lagrangian foliations
([Lib]). A parak\"ahler Lie algebra ${\cal G}$ is the Lie algebra of
a Lie group $G$ with a $G$-invariant parak\"ahler structure ([Ka]).
It is a symplectic Lie algebra with a decomposition into a direct
sum of the underlying vector spaces of two Lagrangian subalgebras.
Some basic facts on the parak\"ahler structures on Lie groups and
Lie algebras have been given in [Bai4]. On the other hand, a phase
space of a Lie algebra in mathematical physics ([Ku1-2], [Bai2]) is
a parak\"ahler Lie algebra. We will show in this paper that
conversely every parak\"ahler Lie algebra is isomorphic to  a phase
space of a Lie algebra.

We have obtained a structure theory of parak\"ahler Lie algebras in
terms of matched pairs of Lie algebras (cf. Theorem 2.5) in [Bai2]
and [Bai4]. This theory in fact gives a construction of parak\"ahler
Lie algebras. However, except for some examples, it is still unclear
when the compatibility conditions appearing in the structure theory
are satisfied.

The aim of this paper is to  study further the structures of
parak\"ahler Lie algebras or phase spaces of Lie algebras in terms
of left-symmetric algebras and interpret the construction mentioned
above using certain equivalent conditions which are much easier to
use. Briefly speaking, a parak\"ahler Lie algebra is equivalent to a
certain bialgebra structure, namely, a left-symmetric bialgebra
structure. From the point of view of phase spaces of Lie algebras,
such a structure seems to be very similar to the Lie bialgebra
structure given by Drinfeld ([D]). In fact, left-symmetric
bialgebras have many properties similar to those of Lie bialgebras.
In particular, there are so-called coboundary left-symmetric
bialgebras which lead to an analogue ($S$-equation) of the classical
Yang-Baxter equation. In a certain sense, the $S$-equation in a
left-symmetric algebra reveals the left-symmetry of the products. A
symmetric solution of the $S$-equation gives a parak\"ahler Lie
algebra.

Furthermore, comparing left-symmetric bialgebras and Lie bialgebras
in terms of several different properties, we observe that there is a
clear analogy between them and in particular, parak\"ahler Lie
groups correspond to Poisson-Lie groups whose Lie algebras are Lie
bialgebras in this sense. Since the classical Yang-Baxter equation
can be regarded as a ``classical limit" of the quantum Yang-Baxter
equation ([Be]), the analogy mentioned above, especially, the
$S$-equation corresponding to the classical Yang-Baxter equation
found in this paper, suggests that there might exist an analogue
(``quantum $S$-equation" ) of the quantum Yang-Baxter equation. The
results in this paper are the beginning of a program to develop the
theory of such analogues of the quantum Yang-Baxter equation. We
expect that our future study will be related to the theory of
quantum groups, tensor categories and vertex operator algebras.

We would like to point out that many structures (for example, see
Theorem 3.8, Theorem 5.4 and so on) appearing in this paper exhibit
features of both Lie algebras and left-symmetric algebras, although
the study of parak\"ahler Lie algebras seems to be purely  a topic
in Lie algebras. Indeed, the theory of Lie algebras alone is not
enough here. Hence, unlike the theory of Lie bialgebras which is
purely Lie-algebra-theoretic, we need to combine the ideas and
methods from both the theory of Lie algebras and the theory of
left-symmetric algebras.

The paper is organized as follows. In Section 2, we give some
necessary definitions and notations and basic results on
left-symmetric algebras and parak\"ahler Lie algebras. In Section 3,
we study how to construct a left-symmetric algebra which is the
direct sum of two left-symmetric subalgebras. We observe that in the
case of parak\"ahler Lie algebras, matched pairs of left-symmetric
algebras are equivalent to the corresponding matched pairs of their
sub-adjacent Lie algebras, whereas it is not true in general. This
also partly explains why left-symmetric algebras appear in a problem
which seems to be purely Lie-algebra-theoretic. In Section 4, we
introduce the notion of left-symmetric bialgebra which is precisely
equivalent to the notion of parak\"ahler Lie algebra. In Section 5,
we consider the special case that a certain 1-cocycle appearing in a
left-symmetric bialgebra is coboundary. A sufficient and necessary
condition for the existence of such a structure leads to certain
explicit equations. In Section 6, we discuss only the simplest cases
in Section 5. We obtain an equation in the left-symmetric algebra,
namely, the $S$-equation, which is an analogue of the classical
Yang-Baxter equation in a Lie algebra. We also give some important
properties of the $S$-equation. In Section 7, we compare
left-symmetric bialgebras and Lie bialgebras by recalling some facts
on Lie bialgebras. We also consider the case that a left-symmetric
bialgebra is also a Lie bialgebra.

Throughout this paper, all algebras are finite-dimensional, although
many results still hold in the infinite-dimensional case.

\section{Preliminaries and basic results}

{\bf Definition 2.1}\quad Let $A$ be a vector space over a field ${\bf F}$
with a bilinear product $(x,y)\rightarrow xy$. $A$ is called a
left-symmetric algebra if for any $x,y,z\in A$, the associator
$$(x,y,z)=(xy)z-x(yz)\eqno (2.1)$$
is symmetric in $x,y$, that is,
$$(x,y,z)=(y,x,z),\;\;{\rm or}\;\;{\rm
equivalently}\;\;(xy)z-x(yz)=(yx)z-y(xz).\eqno (2.2)$$

Left-symmetric algebras are Lie-admissible algebras (cf. [Me]).

{\bf Proposition 2.1}\quad Let $A$ be a left-symmetric algebra. For
any $x,y\in A$, let $L_x$ and $R_x$ denote the left and right
multiplication operator respectively, that is,
$L_x(y)=xy,\;R_x(y)=yx$. Let $L:A\rightarrow gl(A)$ with
$x\rightarrow L_x$ and $R:A\rightarrow gl(A)$ with $x\rightarrow
R_x$ (for every $x\in A$) be two linear maps. Then we have the
following results:

(1) The commutator
$$[x,y]=xy-yx,\;\;\forall x,y\in A,\eqno (2.3)$$
defines a Lie algebra ${\cal G}(A)$, which is called the
sub-adjacent Lie algebra of $A$ and $A$ is also called the
compatible left-symmetric algebra structure on the Lie algebra
${\cal G}(A)$.

(2) $L$ gives a regular representation of the Lie algebra ${\cal G}(A)$, that is,
$$[L_x,L_y]=L_{[x,y]},\;\;\forall x,y\in A. \eqno (2.4)$$

(3) The identity (2.2) is equivalent to the following equation
$$[L_x,R_y]=R_{xy}-R_yR_x,\;\;\forall x,y\in A. \eqno (2.5)$$

Left-symmetric algebras can be obtained from some known algebraic
structures (it can be regarded as a ``realization theory"). Recall
that a Novikov algebra $A$ is a left-symmetric algebra satisfying
$R_xR_y=R_yR_x$ for any $x,y\in A$.

{\bf Example 2.1}\quad  Let $(A,\cdot)$ be a commutative associative
algebra and $D$ be its derivation. Then the new product
$$x*_ay=x\cdot Dy+a\cdot x\cdot y,\;\;\forall\;x,y\in A\eqno (2.6)$$
makes $(A,\;*_a)$ become a Novikov algebra for $a=0$ by S. Gel'fand
([GD]), for $a\in {\bf F}$ by Filipov ([F]) and for a fixed element
$a\in A$ by Xu ([X]). In [BM2-3], we constructed a deformation
theory of Novikov algebras. In particular, the two kinds of Novikov
algebras given by Filipov and Xu are the special deformations of the
algebras $(A,\;*)=(A,\;*_0)$ given by S. Gel'fand. Moreover, we
proved that the Novikov algebras in dimension $\leq 3$ ([BM1]) can
be realized as the algebras defined by S. Gel'fand and their
compatible linear deformations. We conjectured that this conclusion
can be extended to higher dimensions. On the other hand, due to the
structures of free Novikov algebras, any Novikov algebra is a
quotient of a subalgebra of an (infinite-dimensional) algebra given
by equation (2.6) for $a=0$ ([DL]).

{\bf Example 2.2}\quad  Let $({\cal G},[,])$ be a Lie algebra and
$R:{\cal G}\rightarrow {\cal G}$ be a linear map satisfying the
operator form of classical Yang-Baxter equation ([Se] and the Remark
after Proposition 7.5)
$$[R(x),R(y)]=R([R(x),y]+[x,R(y)]),\;\;\forall\; x,y\in {\cal G}.\eqno (2.7) $$ Then
$$x*y=[R(x),y],\;\;\forall\; x,y\in {\cal G}\eqno (2.8)$$
defines a left-symmetric algebra ([GS],[BM5]). Recall that the
matrix form of the above $R$ (satisfying equation (2.7)) is a
classical $r$-matrix. Hence the above construction of left-symmetric
algebras can be regarded as a Lie algebra ``left-twisted" by a
classical $r$-matrix, which gives an algebraic interpretation of
``left-symmetry" (comparing with the geometric interpretation given
in the Introduction) ([Bai3]).

{\bf Example 2.3}\quad Let $(A,\cdot)$ be an associative algebra and
$R:A\rightarrow A$ be a linear map satisfying
$$R(x)\cdot R(y)+R(x\cdot y)=R(R(x)\cdot y+x\cdot R(y)),\;\;\forall
x,y\in A.\eqno (2.9)$$ Then
$$x*y=R(x)\cdot y-y\cdot R(x)-x\cdot y, \;\;\forall
x,y\in A\eqno (2.10)$$ defines a left-symmetric algebra ([GS]). The
above $R$ is called a Rota-Baxter operator which was introduced to
solve analytic ([Bax]) and combinatorial problems ([R]) and attracts
more attention in many fields in mathematics and mathematical
physics ([EGK] and the references therein). It is also related to
the ``modified classical Yang-Baxter equation" ([Se]).

{\bf Example 2.4}\quad Let $V$ be a vector space over the complex
field ${\bf C}$ with the ordinary scalar product $(,)$ and $a$ be a
fixed vector in $V$, then
$$u*v=(u,v)a+(u,a)v,\;\;\forall\; u,v\in V\eqno (2.11)$$
defines a left-symmetric algebra on $V$ which gives the integrable
(generalized) Burgers equation ([SS])
$$U_t=U_{xx}+2U*U_x+(U*(U*U))-((U*U)*U).\eqno (2.12)$$
In [Bai1], we generalized the above construction to get
left-symmetric algebras from linear functions and in particular, we
proved that the left-symmetric algebras given by equation (2.11) are
simple (without any ideals except zero and itself) (cf. [Bu1]).

On the other hand,

{\bf Definition 2.2}\quad A Lie algebra ${\cal G}$ is called a
symplectic Lie algebra if there is a nondegenerate skew-symmetric
2-cocycle $\omega$ (the symplectic form) on ${\cal G}$, that is,
$$\omega([x,y],z)+\omega([y,z],x)+\omega([z,x],y)=0,\;\;\forall x,y,z\in {\cal G}.\eqno (2.13)$$
We denote it by $({\cal G},\omega)$. Let ${\cal H}$ be a subalgebra of a symplectic Lie algebra $({\cal G},\omega)$ and
${\cal H}^\perp =\{ x\in {\cal G}|\omega(x,y)=0,\;\;\forall\; y\in {\cal H}\}$. If ${\cal H}\subset {\cal H}^\perp $, then ${\cal H}$ is called totally isotropic.
If ${\cal H}={\cal H}^\perp $, then ${\cal H}$ is called Lagrangian.

{\bf Theorem 2.2}\quad ([Ch])\quad Let $({\cal G},\omega)$ be a
symplectic Lie algebra. Then there exists a compatible
left-symmetric algebra structure $``*''$ on ${\cal G}$ given by
$$\omega (x*y,z)=-\omega (y,[x,z]),\;\;\forall x,y,z\in {\cal G}.\eqno (2.14)$$

Throughout this paper, we mainly study the following symplectic Lie algebras.

{\bf Definition 2.3}\quad Let $({\cal G},\omega)$ be a symplectic
Lie algebra. ${\cal G}$ is called a parak\"ahler Lie algebra if
${\cal G}$ is a direct sum of the underlying vector spaces of two
Lagrangian subalgebras ${\cal G}^+$ and ${\cal G}^-$. It is denoted
by $({\cal G},{\cal G}^+,{\cal G}^-,\omega)$. Two parak\"ahler Lie
algebras $({\cal G}_1, {\cal G}_1^+,{\cal G}_1^-,\omega_1)$ and
$({\cal G}_2, {\cal G}_2^+,{\cal G}_2^-,\omega_2)$ are isomorphic if
there exists a Lie algebra isomorphism $\varphi:{\cal
G}_1\rightarrow {\cal G}_2$ such that
$$\varphi({\cal G}_1^+)={\cal G}_2^+,\;\;\varphi({\cal G}_1^-)={\cal
G}_2^-,\;\;\omega_1(x,y)=\varphi^*\omega_2(x,y)=\omega_2(\varphi(x),\varphi(y)),\;\;\forall
x,y\in{\cal G}_1.\eqno(2.15)$$

It is easy to see that a symplectic Lie algebra $({\cal G},\omega)$
is a parak\"ahler Lie algebra if and only if ${\cal G}$ is a direct
sum of the underlying vector spaces of two totally isotropic
subalgebras.

{\bf Proposition 2.3}\quad ([Bai2])\quad Let $({\cal G},{\cal
G}^+,{\cal G}^-,\omega)$ be a parak\"ahler Lie algebra, then there
exists a left-symmetric algebra structure on ${\cal G}$ given by
equation (2.14) such that ${\cal G}^+$ and ${\cal G}^-$ are
left-symmetric subalgebras. Moreover, two parak\"ahler Lie algebras
$({\cal G}_i, {\cal G}_i^+,{\cal G}_i^-,\omega_i)$ $(i=1,2)$ are
isomorphic if and only if there exist an isomorphism of
left-symmetric algebras satisfying equation (2.15) which the
compatible left-symmetric algebras are given by equation (2.14).

{\bf Definition 2.4}\quad Let ${\cal G}$ be a Lie algebra. If there
is a Lie algebra structure on a direct sum of the underlying vector
spaces of ${\cal G}$ and its dual space ${\cal G}^*={\rm Hom} ({\cal
G},{\bf F})$ such that ${\cal G}$ and ${\cal G}^*$ are Lie
subalgebras and the natural
 skew-symmetric bilinear form on ${\cal G}\oplus {\cal G}^*$
$$\omega_p (x+a^*,y+b^*)=\langle a^*,y\rangle  -\langle x,b^*\rangle  ,\;\;\forall x,y\in {\cal G},a^*,b^*\in {\cal G}^*,\eqno (2.16)$$
is a 2-cocycle, where $\langle ,\rangle  $ is the ordinary pair
between ${\cal G}$ and ${\cal G}^*$, then it is called a phase space
of the Lie algebra ${\cal G}$ ([Ku1],[Bai2]).

Obviously, a phase space of a Lie algebra is a parak\"ahler Lie algebra. Furthermore, we have

{\bf Proposition 2.4}\quad  Every parak\"ahler Lie algebra $({\cal G},{\cal G}^+,{\cal G}^-,\omega)$
is isomorphic to a phase space of ${\cal G}^+$.

{\bf Proof}\quad Let $\varphi:{\cal G}^-\rightarrow ({\cal G}^+)^*$
be a linear isomorphism given by $\langle \varphi (a),x\rangle
=\omega (a,x)$ for any $x\in {\cal G}^+$ and $a\in {\cal G}^-$.
Extending $\varphi$ to be a linear isomorphism from ${\cal G}={\cal
G}^+\oplus {\cal G}^-$ to ${\cal G}_p={\cal G}^+\oplus ({\cal
G}^+)^*$ with $\varphi(x)=x$ for any $x\in {\cal G}^+$. Let the Lie
algebra structure on ${\cal G}_p$ given by
$$[x,y]=[x,y],\;[a^*,b^*]=\varphi[\varphi^{-1}(a^*),\varphi^{-1}(b^*)],\;\;
[x,a^*]=\varphi([x,\varphi^{-1}(a^*)]),\;\;\forall x\in {\cal G}^+,a^*\in ({\cal G}^+)^*.$$
Therefore $\varphi$ is an isomorphism of Lie algebras. Furthermore, for any $x,y\in {\cal G}^+,a,b\in {\cal G}^-$, we have
$$\omega(x+a,y+b)=\omega(a,y)+\omega(x,b)=\langle \varphi(a),y\rangle  -\langle x,\varphi(b)\rangle
=\omega_p(\varphi(x+a),\varphi(y+b)).$$
Therefore $\omega_p$ is a 2-cocycle in the above Lie algebra ${\cal G}_p$ and
$\varphi$ is an isomorphism of parak\"ahler Lie algebras.
\hfill $\Box$

{\bf Remark}\quad It is obvious that, by symmetry, every parak\"ahler Lie algebra $({\cal G},{\cal G}^+,{\cal G}^-,\omega)$
is also isomorphic to a phase space of ${\cal G}^-$.

On the other hand, recall that $({\cal G},{\cal H}, \rho,\mu)$ is a matched pair of Lie algebras
if ${\cal G}$ and ${\cal H}$ are Lie algebras and $\rho:{\cal G}\rightarrow gl({\cal H})$ and
$\mu:{\cal H}\rightarrow gl({\cal G})$ are representations satisfying
$$\rho(x)[a,b]-[\rho(x)a,b]-[a,\rho(x)b]+\rho(\mu(a)x)b-\rho(\mu(b)x)a=0;\eqno
(2.17)$$
$$\mu(a)[x,y]-[\mu(a)x,y]-[x,\mu(a)y]+\mu(\rho(x)a)y-\mu(\rho(y)a)x=0,\eqno
(2.18)$$ for any $x,y\in {\cal G}$ and $a,b\in {\cal H}$.  In this
case, there exists a Lie algebra structure on the vector space
${\cal G}\oplus {\cal H}$ given by
$$[x+a,y+b]=[x,y]+\mu(a)y-\mu(b)x+[a,b]+\rho(x)b-\rho(y)a,\;\;\forall
x,y\in {\cal G},a,b\in {\cal H}.\eqno (2.19)$$ We denote it by
${\cal G}\bowtie^\rho_\mu {\cal H}$ or simply ${\cal G}\bowtie {\cal
H}$. Moreover, every Lie algebra which is a direct sum of the
underlying vector spaces of two subalgebras can be obtained from a
matched pair of Lie algebras ([Maj], [T]).

{\bf Theorem 2.5}\quad ([Bai2])\quad Let $(A,\cdot)$ be a
left-symmetric algebra. Suppose there is another left-symmetric
algebra structure $``\circ''$ on its dual space $A^*$. Let
$L^*_\cdot$ and $L^*_\circ$ be the dual representation of the
regular representation of ${\cal G}(A)$ and ${\cal G}(A^*)$
respectively, that is,
$$\langle L^*_\cdot (x) a^*,y\rangle  =-\langle a^*,x\cdot y\rangle  ,\;\;\langle L^*_\circ (a^*)x,b^*\rangle  =-\langle x,a^*\circ b^*\rangle  ,\;\;\forall x,y\in A,\;\;
a^*,b^*\in A^*.\eqno (2.20)$$ Then there exists a parak\"ahler Lie
algebra structure on the vector space $A\oplus A^*$ such that ${\cal
G}(A)$ and ${\cal G}(A^*)$ are Lagrangian subalgebras associated to
the symplectic form (2.16) if and only if $({\cal G}(A),{\cal
G}(A^*), L^*_\cdot, L^*_\circ)$ is a matched pair of Lie algebras.
Furthermore, every parak\"ahler Lie algebra can be obtained from the
above way.

\section{Bimodules and matched pairs of left-symmetric algebras}

{\bf Definition 3.1}\quad Let $A$ be a left-symmetric algebra
and $V$ be a vector space. Let $S, T:A\rightarrow gl(V)$ be two linear maps.
$V$ (or the pair $(S,T)$, or $(S,T,V)$) is called a bimodule of $A$ if
$$S(x)S(y)v-S(xy)v=S(y)S(x)v-S(yx)v,\;\;\forall x,y\in A,v\in V;\eqno (3.1)$$
$$S(x)T(y)v-T(y)S(x)v=T(xy)v-T(y)T(x)v,\;\;\forall x,y\in A,v\in V.\eqno (3.2)$$

According to [Sc], the following result is obvious.

{\bf Proposition 3.1}\quad Let $(A,\cdot)$ be a left-symmetric
algebra and $V$ be a vector space. Let $S, T:A\rightarrow gl(V)$ be
two linear maps. Then $(S,T)$ is a bimodule of $A$ if and only if
the direct sum $A\oplus V$ of vector spaces is turned into a
left-symmetric algebra ({the semidirect sum}) by defining
multiplication in $A\oplus V$ by
$$(x_1+v_1)*(x_2+v_2)=x_1\cdot x_2+(S(x_1)v_2+T(x_2)v_1),\;\;\forall x_1,x_2\in, v_1,v_2\in V
\eqno(3.3)$$
We denote it by $A\ltimes_{S,T}V$ or simply $A\ltimes V$.

{\bf Lemma 3.2}\quad Let $(S,T,V)$ be a bimodule of a left-symmetric algebra $A$. Then

(1) $S:A\rightarrow gl(V)$ is a representation of the sub-adjacent Lie algebra ${\cal G}(A)$.

(2) $\rho=S-T$ is a representation of the Lie algebra ${\cal G}(A)$.

(3) For any representation $\mu:{\cal G}(A)\rightarrow gl({\cal G}(A))$ of the Lie algebra ${\cal G}(A)$, $(\mu,0)$ is a bimodule of $A$.

(4) The left-symmetric algebras $A\ltimes_{S,T}V$ and
$A\ltimes_{S-T,0}V$ given by the bimodules $(S,T)$ and $(S-T,0)$
respectively have the same sub-adjacent Lie algebra:
$$[x_1+v_1,x_2+v_2]=[x_1, x_2]+(S-T)(x_1)v_2-(S-T)(x_2)v_1,\;\;\forall x_1,x_2\in {\cal G}(A),v_1,v_2\in V.\eqno (3.4)$$
It is the semidirect sum of a Lie algebra ${\cal G}(A)$ and its module $(\rho=S-T,V)$,
we denote it by ${\cal G}(A)\ltimes_{S-T}V$.

{\bf Proof}\quad (1) and (3) follow directly from Definition 3.1.
For any $x,y\in A$, we have
\begin{eqnarray*}
[(S-T)(x),(S-T)(y)]&=&[S(x),S(y)]-[S(x),T(y)]-[T(x),S(y)]+[T(x),T(y)]\\
&=&S([x,y])-T(xy)+T(y)T(x)+T(yx)-T(x)T(y)+[T(x),T(y)]\\
&=&(S-T)([x,y]).
\end{eqnarray*}
Thus $\rho=S-T$ is a representation of the Lie algebra ${\cal G}(A)$.
Hence (4) follows immediately.
\hfill $\Box$

{\bf Proposition 3.3}\quad Let $(A,\cdot)$ be a left-symmetric algebra and $(S,T,V)$ be its bimodule.
Let $S^*,T^*:A\rightarrow gl(V^*)$ be the linear maps given by
$$\langle S^*(x)u^*,v\rangle  =-\langle S(x)v,u^*\rangle  ,\;\langle T^*(x)u^*,v\rangle  =-\langle T(x)v,u^*\rangle  ,\;\forall x\in A, u^*\in V^*, v\in V.\eqno (3.5)$$
Then $(S^*-T^*,-T^*,V^*)$ is a bimodule of $A$. Therefore there are
two compatible left-symmetric algebra structures
$A\ltimes_{S^*,0}V^*$ and $A\ltimes_{S^*-T^*,-T^*}V^*$ on the Lie
algebra ${\cal G}(A)\ltimes_{S^*}V^*$.

{\bf Proof}\quad Obviously, $S^*-T^*$ is just the dual
representation of the representation $\rho=S-T$ of the sub-adjacent
Lie algebra ${\cal G}(A)$. For any $x,y\in A$, $u^*\in V^*$, $v\in
V$, we have

{\small \begin{eqnarray*}
&&\langle [(S^*-T^*)(x),-T^*(y)]u^*,v\rangle  =\langle (S^*-T^*)(x)u^*,-T(y)v\rangle  +\langle T^*(y)u^*, (S-T)(x)v\rangle  \\
&&=\langle u^*,(S-T)(x)T(y)v\rangle  -\langle u^*,T(y)(S-T)(x)v\rangle  =\langle u^*,(T(xy)-T(x)T(y))v\rangle  \\
&&=\langle -T^*(xy)u^*-T^*(y)T^*(x)u^*,v\rangle  .
\end{eqnarray*}}
Therefore $(S^*-T^*,-T^*,V^*)$  is a bimodule of $A$. The other
conclusions follow from Lemma 3.2 directly.\hfill $\Box$

However, $(S^*,T^*,V^*)$ is not a bimodule of $A$ in general. Similar to the above discussion, we have the following conclusion:

{\bf Proposition 3.4}\quad Let $A$ be a left-symmetric algebra and $(S,T,V)$ be its bimodule. Then the
following conditions are equivalent:

(1) $(S-T,-T)$ is a bimodule of $A$;

(2) $(S^*,T^*)$ is a bimodule of $A$;

(3) $T(x)T(y)=T(y)T(x)$ for any $x,y\in A$.

{\bf Example 3.1}\quad Let $A$ be a left-symmetric algebra. Then $(L,0)$,  $({\rm ad},0)$ and $(L,R)$ are bimodules of $A$, where ${\rm ad}=L-R:A\rightarrow gl(A)$ is the adjoint representation of the Lie algebra ${\cal G}(A)$. On the other hand, $(L^*,0)$, $({\rm ad}^*,0)$ and $({\rm ad}^*,-R^*)$
are bimodules of $A$, too, where ${\rm ad}^*$ and $L^*$ are
the dual representations of the adjoint representation and the regular representation of the sub-adjacent
Lie algebra ${\cal G}(A)$ respectively, and $R^*=-{\rm ad}^*+L^*$.

{\bf Example 3.2}\quad Let $A$ be a left-symmetric algebra. Then
$(L^*,R^*)$ is a bimodule of $A$ if and only if $({\rm ad}, -R)$ is
a bimodule of $A$, if and only if  $A$ is a Novikov algebra.

{\bf Example 3.3}\quad Let $A$ be a left-symmetric algebra. Then
there are two compatible left-symmetric algebras
$A\ltimes_{L^*,0}A^*$ and $A\ltimes_{{\rm ad}^*,-R^*}A^*$ on the Lie
algebra ${\cal G}(A)\ltimes_{L^*}{\cal G}^*(A)$, which are given by
the bimodules $(L^*,0)$ and $({\rm ad}^*,-R^*)$ respectively. This
conclusion has already been given in [Bai2].

Next, we study how to construct a left-symmetric algebra structure
on a direct sum $A\oplus B$ of the underlying vector spaces of two
left-symmetric algebras $A$ and $B$ such that $A$ and $B$ are
subalgebras.

{\bf Theorem 3.5}\quad  Let $(A,\cdot)$ and $(B,\circ)$ be two left-symmetric algebras.
Suppose that there are linear maps
$l_A,r_A:A\rightarrow gl(B)$ and $l_B,r_B:B\rightarrow gl(A)$ such that
$(l_A,r_A)$ is a bimodule of $A$ and $(l_B,r_B)$ is a bimodule of $B$ and they satisfy the following conditions:
$$r_A(x)[a,b]=r_A(l_B(b)x)a-r_A(l_B(a)x)b+a\circ (r_A(x)b)-b\circ (r_A(x)a);\eqno (3.6)$$
$$l_A(x)(a\circ b)=-l_A(l_B(a)x-r_B(a)x)b+(l_A(x)a-r_A(x)a)\circ b+r_A(r_B(b)x)a+a\circ(l_A(x)b);\eqno (3.7)$$
$$r_B(a)[x,y]=r_B(l_A(y)a)x-r_B(l_A(x)a)y+x\cdot (r_B(a)y)-y\cdot (r_B(a)x);\eqno (3.8)$$
$$l_B(a)(x\cdot y)=-l_B(l_A(x)a-r_A(x)a)y+(l_B(a)x-r_B(a)x)\cdot y+r_B(r_A(y)a)x+x\cdot(l_B(a)y),\eqno (3.9)$$
for any $x,y\in A,a,b\in B$. Then there is a left-symmetric algebra
structure on the vector space $A\oplus B$ given by
$$(x+a)*(y+b)=(x\cdot y+l_B(a)y+r_B(b)x)+(a\circ b+l_A(x)b+r_A(y)a),\;\;\forall x,y\in A,a,b\in B.\eqno (3.10)$$
We denote this left-symmetric algebra by
$A\bowtie^{l_A,r_A}_{l_B,r_B}B$ or simply $A\bowtie B$. And
$(A,B,l_A,r_A,l_B,r_B)$ satisfying the above conditions is called a
matched pair of left-symmetric algebras. On the other hand, every
left-symmetric algebra which is a direct sum of the underlying
vector spaces of two subalgebras can be obtained from the above way.

{\bf Proof}\quad For any $x,y,z\in A$ and $a,b,c\in B$, the associator in $A\bowtie B$ satisfies
$$(x+a,y+b,z+c)=(x,y,z)+(x,b,c)+(a,y,c)+(x,y,c)+(x,b,z)+(a,b,c)+(a,y,z)+(a,b,z).$$
Therefore, the equation (3.10) defines a left-symmetric algebra
structure on $A\oplus B$ if and only if the following equations are
satisfied:
\begin{eqnarray*}
 (x,y,z)=(y,x,z)\;&\Longleftrightarrow&\;\;A\;{\rm is}\;{\rm a}\;{\rm left-symmetric}\;{\rm algebra};\\
 (x,y,c)=(y,x,c)\;&\Longleftrightarrow&\;\;l_A\;\; {\rm is}\;{\rm a}\;
{\rm representation}\;{\rm of}\;{\cal G}(A)\;{\rm and}\;{\rm equation}\; (3.8)\;{\rm holds};\\
(a,y,z)=(y,a,z)\;&\Longleftrightarrow&\;\; [(x,b,z)=(b,x,z)\;\Longleftrightarrow] \;\;{\rm equation}\; (3.9)\;{\rm holds}\;{\rm and}\;\;r_A\;\; {\rm satisfies}\\
&\mbox{}&\;\; l_A(y)r_A(z)a-r_A(z)l_A(y)a=r_A(y\cdot z)a-r_A(z)r_A(y)a;\\ (a,b,c)=(b,a,c)\;&\Longleftrightarrow&\;\;B\;{\rm is}\;{\rm a}\;{\rm left-symmetric}\;{\rm algebra};\\
(a,b,z)=(b,a,z)\;&\Longleftrightarrow&\;\;l_B\;\; {\rm is}\;{\rm a}\;
{\rm representation}\;{\rm of}\;{\cal G}(B)\;{\rm and}\;{\rm equation}\; (3.6)\;{\rm holds};\\
(x,b,c)=(b,x,c)\;&\Longleftrightarrow&\;\; [(a,y,c)=(y,a,c)\;\Longleftrightarrow] \;\;{\rm equation}\; (3.7)\;{\rm holds}\;{\rm and}\;\;r_B\;\; {\rm satisfies}\\
&\mbox{}& l_B(b)r_B(c)x-r_B(c)l_B(b)x=r_B(b\circ c)x-r_B(c)r_B(b)x.
\end{eqnarray*}
Hence $A\bowtie B$ is a left-symmetric algebra if and only if
$(l_A,r_A)$ is a bimodule of $A$ and $(l_B,r_B)$ is a bimodule of
$B$ and equations (3.6-3.9) hold. On the other hand, if $A$ and $B$
are left-symmetric subalgebras of a left-symmetric $C$ such that
$C=A\oplus B$ which is a direct sum of the underlying vector spaces
of $A$ and $B$, then it is easy to show that the linear maps
$l_A,r_A:A\rightarrow gl(B)$ and $l_B,r_B:B\rightarrow gl(B)$
determined by
$$x*a=l_A(x)a+r_B(a)x,\;\;a*x=l_B(a)x+r_A(x)a,\;\;\forall x\in A,a\in B$$
satisfy equations (3.6-3.9). In addition, $(l_A,r_A)$ is a bimodule of $A$ and $(l_B,r_B)$ is a
bimodule of $B$.\hfill $\Box$

{\bf Remark}\quad Obviously $B$ is an ideal of $A\bowtie B$ if and only if $l_B=r_B=0$. In this case,
if in addition, $B$ is trivial (that is, all the products of $B$ are zero),
then $A\bowtie^{l_A,r_A}_{0,0} B\cong A\ltimes_{l_A,r_A} B$. Furthermore,
some other special cases of the above theorem have already been studied. For example, Diatta and Medina studied
the case that both $A$ and $B$ are left ideals of $A\bowtie B$, that is, $r_A=r_B=0$ ([DiM]).

{\bf Corollary 3.6} ([DiM, Lemma 3.0.3])\quad Let $(A,\cdot)$ and $(B,\circ)$ be two left-symmetric algebras. Suppose there
exist two Lie algebra representations $\rho:{\cal G}(A)\rightarrow gl(B)$ and
$\mu:{\cal G}(B)\rightarrow gl(A)$ such that
$$\rho(x)(a\circ b)=\rho(x)a\circ b+a\circ (\rho(x)b)-\rho(\mu(a)x)b;\eqno (3.11)$$
$$\mu(a)(x\cdot y)=\mu(a)x\cdot y+x\cdot (\mu(a)y)-\mu(\rho(x)a)y,\eqno (3.12)$$
for any $x,y \in A,a,b\in B$. $A$ and $B$ are called to be
$(\rho,\mu)$-linked. Then there is a left-symmetric algebra
structure on the vector space $A\oplus B$ given by
$$(x+a)*(y+b)=(x\cdot y+\mu (a)y)+(a\circ b+\rho(x)b),\forall x,y\in A,a,b\in B.\eqno (3.13)$$
Therefore, its sub-adjacent Lie algebra (which is the same Lie
algebra given by the matched pair of Lie algebras $({\cal
G}(A),{\cal G}(B),\rho,\mu)$) is given by
$$[(x+a),(y+b)]=(x\cdot y-y\cdot x+\mu (a)y-\mu(b)x)+(a\circ b-b\circ a +\rho(x)b-\rho(y)a),\forall x,y\in A,a,b\in B.\eqno (3.14)$$

{\bf Corollary 3.7}\quad Let $(A,B,l_A,r_A,l_B,r_B)$ be a matched pair of left-symmetric algebras. Then
$({\cal G}(A),{\cal G}(B),l_A-r_A,l_B-r_B)$ is a matched pair of Lie algebras.

{\bf Proof}\quad This conclusion can be proved by a direct
computation or from the relation between the left-symmetric algebra
$A\bowtie B$ and its sub-adjacent algebra Lie algebra. In fact, the
sub-adjacent Lie algebra ${\cal G}(A\bowtie B)$ is just the Lie
algebra obtained from the matched pair $({\cal G}(A),{\cal
G}(B),\rho,\mu)$:
$$[x+a,y+b]=[x,y]+\mu(a)y-\mu(b)x+[a,b]+\rho(x)b-\rho(y)a,\;\;\forall
x,y\in {\cal G}(A),a,b\in {\cal G}(B),$$
where $\rho=l_A-r_A,\mu=l_B-r_B$.\hfill $\Box$

Now we turns to the case of a parak\"ahler Lie algebra: $\rho=L^*_\cdot$ and $\mu=L^*_\circ$.

{\bf Theorem 3.8}\quad Let $(A,\cdot)$ be a left-symmetric algebra.
Suppose there is another left-symmetric algebra structure
$``\circ''$ on its dual space $A^*$. Then $({\cal G}(A),{\cal
G}(A^*),L^*_\cdot,L^*_\circ)$ is a matched pair of Lie algebras if
and only if $(A,A^*, {\rm ad}^*_\cdot, -R^*_\cdot,{\rm
ad}^*_\circ,-R^*_\circ)$ is a matched pair of left-symmetric
algebras.

{\bf Proof }\quad We only need to prove the ``only if" part. A direct proof is given as follows. In fact,
in the case $l_A={\rm ad}^*_\cdot, r_A=-R^*_\cdot$, $l_B=l_{A^*}={\rm ad}^*_{\circ},r_B=r_{A^*}=-R^*_\circ$,
$\rho=L^*_\cdot$ and $\mu=L^*_\circ$, we have
$${\rm equation}\;\; (2.17) \;\; \Longleftrightarrow {\rm equation}\;\; (3.6)\;\;
\Longleftrightarrow\;\; {\rm equation}\;\; (3.9);$$
$${\rm equation}\;\; (2.18) \;\; \Longleftrightarrow {\rm equation}\;\; (3.7)\;\;
\Longleftrightarrow\;\; {\rm equation}\;\; (3.8).$$ As an example
(the proof of other equivalent relations is similar), we show how
equation (2.17) is equivalent to equation (3.6). In fact, it follows
from
\begin{eqnarray*}
&&\langle -R^*_\cdot(x)[a^*,b^*],y\rangle  =-\langle L^*(y)[a^*,b^*],x\rangle  ;\\
&&\langle -R^*_\cdot({\rm ad}^*_{\circ}(b^*)x)a^*,y\rangle  =\langle [b^*,L^*_\cdot(y)a^*],x\rangle  ;\\
&&-\langle -R^*_\cdot({\rm ad}^*_{\circ}(a^*)x)b^*,y\rangle  =\langle -[a^*,L^*_\cdot(y)b^*],x\rangle  ;\\
&&\langle a^*\circ(-R^*_\cdot(x)b^*),y\rangle  =\langle L^*_\cdot(L^*_\circ(a^*)y)(b^*),x\rangle  ;\\
&&\langle -b^*\circ(-R^*_\cdot(x)a^*),y\rangle  =\langle
-L^*_\cdot(L^*_\circ(b^*)y)(a^*),x\rangle  ,
\end{eqnarray*}
for any $x,y\in A$ and $a^*,b^*\in A^*$.

There is another proof. If $({\cal G}(A),{\cal
G}(A^*),L^*_\cdot,L^*_\circ)$ is a matched pair of Lie algebras,
then ${\cal G}(A)\bowtie^{L^*_\cdot}_{L^*_\circ} {\cal G}(A^*)$ is a
parak\"ahler Lie algebra with the natural symplectic form (2.16).
Hence there exists a compatible left-symmetric algebra structure on
${\cal G}(A)\bowtie^{L^*_\cdot}_{L^*_\circ} {\cal G}(A^*)$ given by
equation (2.14). With a simple and direct computation, we know that
$A$ and $A^*$ are its subalgebras and the other products are given
by
$$x*a^*={\rm ad}^*_\cdot(x) a^*-R^*_\circ(a^*)(x),\;\;a^**x={\rm ad}^*_\circ(a^*)x-R^*_\cdot(x)a^*,\;\forall x\in A,a^*
\in A^*.$$
Therefore $(A,A^*, {\rm ad}^*_\cdot, -R^*_\cdot,{\rm ad}^*_\circ,-R^*_\circ)$
is a matched pair of left-symmetric algebras. \hfill $\Box$

{\bf Remark}\quad Unlike the case of semidirect sum in Example 3.3,
it is not true that there is a compatible left-symmetric structure
on the Lie algebra ${\cal G}(A)\bowtie^{L^*_\cdot}_{L^*_\circ} {\cal
G}(A^*)$ coming from the bimodules $(L^*_\cdot,0)$ and
$(L^*_\circ,0)$ since in this case,
$(A,A^*,L^*_\cdot,0,L^*_\circ,0)$ is not a matched pair of
left-symmetric algebras in general. In fact, from Corollary 3.6,
$(A,A^*,L^*_\cdot,0,L^*_\circ,0)$ is a matched pair of
left-symmetric algebras if and only if
$$L^*_\cdot (x)(a^*\circ b^*)=-L^*_\cdot(L^*_\circ(a^*)x)b^*+L^*_\cdot(x)a^*\circ b^*+a^*\circ L^*_\cdot(x)b^*;\eqno(3.15)$$
$$L^*_\circ(a^*)(x\cdot y)=-L^*_\circ (L^*_\cdot(x)a^*)y+L^*_\circ(a^*)x\cdot y+x\cdot L^*_\circ(a^*)y,\eqno(3.16) $$
for any $x,y\in A$ and $a^*,b^*\in A^*$ which is a little stronger
than equations (2.17-2.18) in the case $\rho=L_\cdot^*$ and
$\mu=L_\circ^*$.

\section{Left-symmetric bialgebras}

At first, we give some notations as follows.

Let ${\cal G}$ be a Lie algebra. For any two representations $(\rho,V)$ and $(\mu,W)$ of ${\cal G}$, it is easy to
know that $(\rho\otimes 1+1\otimes \mu,V\otimes W)$ is also a representation of ${\cal G}$, where
$$(\rho\otimes 1+1\otimes \mu)(x) (v\otimes w)=(\rho(x)\otimes 1+1\otimes \mu(x))(v\otimes w)=
\rho(x)v\otimes w+v\otimes \mu(x) w,\eqno (4.1)$$
for any $x\in {\cal G},v\in V,w\in W$.

For a Lie algebra ${\cal G}$ and a representation $(\rho, V)$ of ${\cal G}$, recall that a 1-cocycle
$\delta$ associated to $\rho$ (denoted by $(\rho, \delta)$) is a linear map from ${\cal G}$ to $V$  satisfying
$$\delta([x,y])=\rho(x)\delta(y)-\rho(y)\delta(x),\;\;\forall x,y\in {\cal G}.\eqno (4.2)$$

For a linear map $\phi:V_1\rightarrow V_2$, we denote the dual (linear) map by $\phi^*:V_2^*\rightarrow V_1^*$ given by
$$\langle v,\phi^*(u^*)\rangle  =\langle \phi(v),u^*\rangle  ,\;\;\forall v\in V_1,u^*\in V_2.\eqno (4.3)$$

{\bf Theorem 4.1}\quad Let $(A,\cdot)$ be a left-symmetric algebra
whose product is given by a linear map $\beta^*:A\otimes
A\rightarrow A$. Suppose there is another left-symmetric algebra
structure $``\circ''$ on its dual space $A^*$ given by a linear map
$\alpha^*: A^*\otimes A^*\rightarrow A^*$. Then $({\cal G}(A),{\cal
G}(A^*), L^*_\cdot, L^*_\circ)$ is a matched pair of Lie algebras if
and only if $\alpha:A\rightarrow A\otimes A $ is a 1-cocycle of
${\cal G}(A)$ associated to $L_\cdot \otimes 1+1\otimes {\rm
ad}_\cdot$ and $\beta: A^*\rightarrow A^*\otimes A^*$ is a 1-cocycle
of ${\cal G}(A^*)$ associated to $L_\circ\otimes 1 + 1\otimes {\rm
ad}_\circ$.

{\bf Proof}\quad Let $\{e_1,\cdots,e_n\}$ be a basis of $A$ and $\{
e_1^*,\cdots, e_n^*\}$ be its dual basis. Set $e_i\cdot
e_j=\sum\limits_{k=1}^nc_{ij}^ke_k$ and $e_i^*\circ
e_j^*=\sum\limits_{k=1}^nf_{ij}^k e_k^*$. Therefore, we have
$$\alpha (e_k)=\sum_{i,j=1}^nf_{ij}^ke_i\otimes e_j,\;\;\beta(e_k^*)=\sum_{i,j=1}^nc_{ij}^k e_i^*\otimes e_j^*;\;\;
L^*_\cdot (e_i)e_j^*=-\sum_{k=1}^nc_{ik}^je_k^*,\;\;L^*_\circ (e_i^*)e_j=-\sum_{k=1}^nf_{ik}^je_k.$$
Hence the coefficient of $e_m\otimes e_n$ in
$$\alpha([e_i,e_j])=[L_\cdot (e_i)\otimes 1+1\otimes {\rm ad}_\cdot (e_i)]\alpha(e_j)-[L_\cdot ({e_j})\otimes 1+
1\otimes {\rm ad}_\cdot (e_j)]\alpha(e_i)$$
gives the following relation (for any $i,j,m,n$)
$$\sum_{k=1}^n(c_{ij}^k-c_{ji}^k)f_{mn}^k=\sum_{k=1}^n[
c_{ik}^mf_{kn}^j+(c_{ik}^n-c_{ki}^n)f_{mk}^j-c_{jk}^mf_{kn}^i-(c_{jk}^n-c_{kj}^n)f_{mk}^i],$$
which is precisely the relation given by the coefficient of $e_n$ in
$$-L^*_\circ(e_m^*)[e_i,e_j]= L^*_\circ (L^*_\cdot(e_i)e_m^*)e_j-
[e_i,L^*_\circ(e_m^*)e_j]-L^*_\circ
(L^*_\cdot(e_j)e_m^*)e_i-[L^*_\circ(e_m^*)e_i,e_j].$$ Then $\alpha$
is a 1-cocycle of ${\cal G}(A)$ associated to $L_\cdot \otimes
1+1\otimes {\rm ad}_\cdot$ if and only if equation (2.18) holds in
the case $\rho=L^*_\cdot$ and $\mu=L^*_\circ$. Similarly (or by
symmetry), $\beta$ is a 1-cocycle of ${\cal G}(A^*)$ associated to
$L_\circ\otimes 1 + 1\otimes {\rm ad}_\circ$ if and only if equation
(2.17) holds.\hfill $\Box$

{\bf Definition 4.1}\quad  Let $A$ be a vector space. A
left-symmetric bialgebra structure on $A$ is a pair of linear maps
$(\alpha,\beta)$ such that $\alpha:A\rightarrow A\otimes A,\beta:
A^*\rightarrow A^*\otimes A^*$ and

(a) $\alpha^*:A^*\otimes A^*\rightarrow A^*$ is a left-symmetric
algebra structure on $A^*$;

(b) $\beta^*:A\otimes A\rightarrow A$ is a left-symmetric algebra
structure on $A$;

(c) $\alpha$ is a 1-cocycle of ${\cal G}(A)$ associated to $L\otimes 1+1\otimes {\rm ad}$ with values in $A\otimes A$;

(d) $\beta$ is a 1-cocycle of ${\cal G}(A^*)$ associated to
$L\otimes 1+1\otimes {\rm ad}$ with values in $A^*\otimes A^*$.

We also denote this left-symmetric bialgebras by $(A,A^*,\alpha,\beta)$ or simply $(A,A^*)$.

{\bf Proposition 4.2}\quad  Let $(A,\cdot)$ be a left-symmetric
algebra and $(A^*,\circ)$ be a left-symmetric algebra structure on
its dual space $A^*$. Then the following conditions are equivalent:

(1) $({\cal G}(A)\bowtie {\cal G}(A)^*, {\cal G}(A),{\cal
G}(A^*),\omega_p)$ is a parak\"ahler Lie algebra, where $\omega_p$
is given by equation (2.16);

(2)  $({\cal G}(A),{\cal G}(A^*),L^*_\cdot,L^*_\circ)$ is
a matched pair of Lie algebras;

(3) $(A,A^*, {\rm ad}^*_\cdot, -R^*_\cdot,{\rm ad}^*_\circ,-R^*_\circ)$ is
a matched pair of left-symmetric algebras;

(4) $(A,A^*)$ is a left-symmetric bialgebra.

{\bf Proof}\quad The proofs are straightforward. \hfill $\Box$

{\bf Definition 4.2}\quad Let $(A,A^*,\alpha_A,\beta_A)$ and $(B,B^*,\alpha_B,\beta_B)$ be
two left-symmetric bialgebras. A homomorphism of left-symmetric bialgebras $\varphi:A\rightarrow B$
is a homomorphism of left-symmetric algebras such that
$\varphi^*:B^*\rightarrow A^*$ is also a homomorphism of left-symmetric algebras, that is, $\varphi$ satisfies
$$(\varphi\otimes \varphi) \alpha_A (x)=\alpha_B(\varphi(x)),\;\;
(\varphi^*\otimes
\varphi^*)\beta_B(a^*)=\beta_A(\varphi^*(a^*)),\forall\; x\in
A,a^*\in B^*.\eqno (4.4)$$ An isomorphism of left-symmetric
bialgebras is an invertible homomorphism of left-symmetric
bialgebras.

{\bf Proposition 4.3}\quad Two parak\"ahler Lie algebras are isomorphic if and only if their corresponding left-symmetric bialgebras
are isomorphic.

{\bf Proof}\quad Let $({\cal G}(A)\bowtie{\cal G}(A^*), {\cal G}(A), {\cal G}(A^*),\omega_p)$ and
$({\cal G}(B)\bowtie{\cal G}(B^*), {\cal G}(B), {\cal G}(B^*),\omega_p)$ be two parak\"ahler Lie algebras.
Let $\{e_1,\cdots,e_n\}$ be a basis of $A$ and $\{e_1^*,\cdots, e_n^*\}$ be its dual basis. If $\varphi:{\cal G}(A)\bowtie {\cal G}(A^*)\rightarrow {\cal G}(B)\bowtie {\cal G}(B^*)$ is an isomorphism of parak\"ahler Lie algebras, then $\varphi|_A:A\rightarrow B$ and $\varphi|_{A^*}:A^*\rightarrow B^*$ are isomorphisms of left-symmetric algebras due to Proposition 2.3. Moreover,
$\varphi|_{A^*}={(\varphi|_{A})^*}^{-1}$ since
\begin{eqnarray*}
\langle \varphi|_{A^*}(e_i^*), \varphi (e_j)\rangle
&=&\omega_p(\varphi|_{A^*}(e_i^*),\varphi(e_j))
=\omega_p(e_i^*,e_j)\\
&=&\delta_{ij}=\langle e_i^*,e_j\rangle  =\langle \varphi^*{(\varphi|_A)^*}^{-1}(e_i^*),e_j\rangle  \\
&=&\langle {(\varphi|_A)^*}^{-1}(e_i^*),\varphi(e_j)\rangle  .
\end{eqnarray*}
Hence $(A,A^*)$ and $(B,B^*)$ are isomorphic as left-symmetric bialgebras.
On the other hand, let $(A,A^*)$ and $(B,B^*)$ be two isomorphic left-symmetric bialgebras and $\varphi':A\rightarrow B$ be an isomorphism of left-symmetric bialgebras. Set $\varphi:A\oplus A^*\rightarrow B\oplus B^*$ be a linear map given by
$$\varphi(x)=\varphi'(x),\varphi(a^*)=(\varphi'^*)^{-1}(a^*),\;\;\forall x\in A,a^*\in A^*.$$
Then it is easy to know that $\varphi$ is an isomorphism of the two
parak\"ahler Lie algebras $({\cal G}(A)\bowtie{\cal G}(A^*), {\cal
G}(A), {\cal G}(A^*),\omega_p)$ and $({\cal G}(B)\bowtie{\cal
G}(B^*), {\cal G}(B), {\cal G}(B^*),\omega_p)$.\hfill $\Box$

{\bf Example 4.1}\quad Let $(A,A^*,\alpha,\beta)$ be a left-symmetric bialgebra. Then its dual $(A^*,A,\beta,\alpha)$ is also a left-symmetric bialgebra.

{\bf Example 4.2}\quad Let $A$ be a left-symmetric algebra. Let the
left-symmetric algebra structure on $A^*$ be trivial, then in this
case $(A,A^*,0,\beta)$ is a left-symmetric bialgebra. And its
corresponding left-symmetric algebra is $A\ltimes_{{\rm
ad}^*,-R^*}A^*$. Moreover, its corresponding  parak\"ahler Lie
algebra is just the semidirect sum ${\cal G}(A)\ltimes_{L^*} {\cal
G}(A^*)$ with the symplectic form $\omega_p$ given by equation
(2.16). Dually, let $A$ be a trivial left-symmetric algebra, then
the left-symmetric bialgebra structures on $A$ are in one-to-one
correspondence with the left-symmetric algebra structure on $A^*$.

{\bf Example 4.3}\quad Let $(A,A^*)$ be a left-symmetric bialgebra.
In next section, we will prove that there exists a natural
left-symmetric bialgebra structure on the direct sum $A\oplus A^*$
of the underlying vector spaces of $A$ and $A^*$.

\section{ Coboundary left-symmetric bialgebras}

In this section, we study the case that $\alpha$ is a 1-coboundary associated to
$L\otimes 1+1\otimes {\rm ad}$.

{\bf Definition 5.1}\quad A left-symmetric bialgebra
$(A,A^*,\alpha,\beta)$ is called coboundary if $\alpha$ is a 1-coboundary of ${\cal G}(A)$
associated to $L\otimes 1+1\otimes {\rm ad}$, that is, there exists a $r\in A\otimes A$ such that
$$\alpha (x)=(L_x\otimes 1+1\otimes {\rm ad}x)r,\;\;\forall x\in A.\eqno (5.1)$$

Let $A$ be a left-symmetric algebra whose product is given by $\beta^*:A\otimes A\rightarrow A$
and $r\in A\otimes A$. Suppose $\alpha:A\rightarrow A\otimes A$ is a 1-coboundary of ${\cal G}(A)$
associated to $L\otimes 1+1\otimes {\rm ad}$ (given by equation (5.1)). Then it is obvious that $\alpha$ is a 1-cocycle of ${\cal G}(A)$ associated to $L\otimes 1+1\otimes {\rm ad}$. Therefore, $(A,A^*,\alpha,\beta)$ is a left-symmetric bialgebra if and only if the
following two conditions are satisfied:

(1) $\alpha^*:A^*\otimes A^*\rightarrow A^*$ defines a
left-symmetric algebra structure on $A^*$.

(2) $\beta$ is a 1-cocycle of ${\cal G}(A^*)$ associated to $L\otimes 1+1\otimes {\rm ad}$, where the left-symmetric algebra
structure of $A^*$ is given by (1).

For any vector space $A$,
let $\sigma:A\otimes A\rightarrow A\otimes A$ be the linear map (exchanging operator) satisfying
$\sigma (x\otimes y)=y\otimes x$, for any $x,y\in A$. Let $r=\sum_{i}a_i\otimes b_i\in A\otimes A$, where $a_i,b_i\in A$.
We denote $r_{12}=r$ and $r_{21}=\sigma (r)=\sum_{i} b_i\otimes a_i$.

{\bf Proposition 5.1}\quad Let $(A,\cdot)$ be a left-symmetric
algebra whose product is given by $\beta^*:A\otimes A\rightarrow A$
and $r\in A\otimes A$. Suppose there exists a left-symmetric algebra
structure $``\circ"$ on $A^*$ given by $\alpha^*:A^*\otimes
A^*\rightarrow A^*$, where $\alpha$ is given by equation (5.1). Then
$\beta:A^*\rightarrow A^*\otimes A^*$ is a 1-cocycle of ${\cal
G}(A^*)$ associated to $L\otimes 1+1\otimes {\rm ad}$ if and only if
$r$ satisfies
$$[P(x\cdot y)-P(x)P(y)](r_{12}-r_{21})=0,\forall x,y\in A,\eqno (5.2)$$
where $P(x)=L_x\otimes 1+1\otimes L_x$.

{\bf Proof}\quad Let $\{ e_1,\cdots, e_n\}$ be a basis of $A$ and $\{ e_1^*,\cdots,e_n^*\}$ be its dual basis.
Set $r=\sum\limits_{i,j} a_{ij}e_i\otimes e_j$. Suppose $e_i\cdot e_j=\sum_{k}c_{ij}^ke_k$ and $e_i^*\circ e_j^*=\sum_{k}f_{ij}^ke_k^*$. Since
$$\alpha (e_i)=\sum_{k,l} f_{kl}^i e_k\otimes e_l=(L_{e_i}\otimes 1+1\otimes {\rm ad}e_i)r,$$
we have (for any $i,k,l$)
$$f_{kl}^i=\sum_{t}[a_{tl}c_{it}^k+a_{kt}(c_{it}^l-c_{ti}^l)].$$
Therefore $\beta$ is a 1-cocycle of ${\cal G}(A^*)$ associated to
$L_\circ\otimes 1+1\otimes {\rm ad}_\circ$ if and only if
$$\beta[e_i^*,e_j^*]=(L_{e_i^*}\otimes 1+1\otimes {\rm ad}e_i^*)\beta^*(e_j^*)-
(L_{e_j^*}\otimes 1+1\otimes {\rm ad}e_j^*)\beta^*(e_i^*),$$ where
$\beta(e_i^*)=\sum_{k,l}c_{kl}^ie_k^*\otimes e_l^*$. Let both sides
of the above equation act on $e_m\otimes e_n$ and after rearranging
the terms suitably, we have
$$(F1)+(F2)+(F3)+(F4)+(F5)+(F6)=0,$$
where
\begin{eqnarray*}
(F1)&=&\sum_{t,l}(a_{tl}-a_{lt})(c_{nt}^jc_{ml}^i-c_{nt}^ic_{ml}^j);\\
(F2)&=&\sum_{t,l}(-a_{lt}c_{tn}^ic_{ml}^j+a_{tl}c_{mt}^jc_{ln}^i-a_{tl}c_{mt}^ic_{ln}^j+a_{lt}c_{ml}^ic_{tn}^j);\\
(F3)&=&\sum_{t,l}(a_{tj}-a_{jt})(c_{lt}^ic_{mn}^l-c_{nt}^lc_{ml}^i);\\
(F4)&=&\sum_{t,l} a_{jt}[c_{tl}^ic_{mn}^l-c_{tn}^lc_{ml}^i+(c_{mt}^l-c_{tm}^l)c_{ln}^i];\\
(F5)&=&\sum_{t,l}(a_{ti}-a_{it})(c_{nt}^lc_{ml}^j-c_{lt}^jc_{mn}^l);\\
(F6)&=&\sum_{t,l} a_{it}[-c_{tl}^jc_{mn}^l+c_{tn}^lc_{ml}^j-(c_{mt}^l-c_{tm}^l)c_{ln}^j].\\
\end{eqnarray*}
$(F1)$ is the coefficient of $e_i\otimes e_j$ in
$$(-L_{e_m}\otimes L_{e_n}-L_{e_n}\otimes L_{e_m})\sum_{t,l}(a_{tl}-a_{lt})e_t\otimes e_l=
(-L_{e_m}\otimes L_{e_n}-L_{e_n}\otimes L_{e_m}](r_{12}-r_{21});$$
$(F2)=0$ by interchanging the indices $t$ and $l$;

\noindent $(F3)$ is the coefficient of $e_i\otimes e_j$ in
$$(L_{e_m\cdot e_n}\otimes 1- L_{e_m}L_{e_n}\otimes 1)\sum_{t,l}(a_{tl}-a_{lt})e_t\otimes e_l=
[L_{e_m\cdot e_n}\otimes 1-L_{e_m}L_{e_n}\otimes
1](r_{12}-r_{21});$$ $(F4)=0$ since the term in the bracket is the
coefficient of $e_i$ in $(e_m,e_t,e_n)-(e_t,e_m,e_n)=0$. $(F5)$ is
the coefficient of $e_i\otimes e_j$ in
$$(1\otimes L_{e_m\cdot e_n}- 1\otimes L_{e_m}L_{e_n})\sum_{t,l}(a_{tl}-a_{lt})e_t\otimes e_l=
[1\otimes L_{e_m\cdot e_n}-1\otimes
L_{e_m}L_{e_n}](r_{12}-r_{21});$$ $F(6)=0$ since the term in the
bracket is the coefficient of $e_j$ in
$(e_t,e_m,e_n)-(e_m,e_t,e_n)=0$. Therefore we have
$$[P(e_m\cdot e_n)-P(e_m)P(e_n)](r_{12}-r_{21})=0.$$
Hence the conclusion holds.\hfill$\Box$

For any linear map $\alpha:A\rightarrow A\otimes A$, let $J_\alpha:A\rightarrow A\otimes A\otimes A$ be a linear map given by
$$J_{\alpha}(x)=(\alpha\otimes {\rm id})\alpha(x)-({\rm id}\otimes \alpha)\alpha(x)-(\sigma\otimes {\rm id})(\alpha\otimes {\rm id})\alpha(x)+(\sigma\otimes {\rm id})({\rm id}\otimes \alpha)\alpha(x),\forall x\in A.\eqno (5.3)$$

{\bf Lemma 5.2}\quad Let $A$ be a vector space and $\alpha:A\otimes
A\rightarrow A$ be a linear map. Then $\alpha^*:A^*\otimes
A^*\rightarrow A^*$ defines a left-symmetric algebra structure on
$A^*$ if and only if $J_\alpha=0$.

{\bf Proof}\quad We denote the product in $A^*$ by $``\circ''$. Then for any $a^*,b^*\in A^*$, we have
$$a^*\circ b^*=\alpha^*(a^*\otimes b^*).$$
Hence the associator satisfies (for any $a^*,b^*,c^*\in A^*$)
\begin{eqnarray*}
(a^*,b^*,c^*)&=&[\alpha^*(\alpha^*\otimes {\rm id})-\alpha^*({\rm id}\otimes
\alpha^*)](a^*\otimes b^*\otimes c^*);\\
(b^*,a^*,c^*)&=&[\alpha^*(\alpha^*\otimes {\rm id})-\alpha^*({\rm id}\otimes
\alpha^*)](b^*\otimes a^*\otimes c^*)\\
&=&[\alpha^*(\alpha^*\otimes {\rm id})-\alpha^*({\rm id}\otimes
\alpha^*)](\sigma\otimes {\rm id})(a^*\otimes b^*\otimes c^*).
\end{eqnarray*}
Therefore $(a^*,b^*,c^*)=(b^*,a^*,c^*)$ for any $a^*,b^*,c^*\in A^*$
if and only if $J_\alpha=0$.\hfill $\Box$

{\bf Proposition 5.3}\quad Let $(A,\cdot)$ be a left-symmetric algebra. Let $r=\sum_{i}a_i\otimes b_i\in A\otimes A$, where
$a_i,b_i\in A$. Define $\alpha:A\rightarrow A\otimes A$ by equation (5.1). Then
$$J_\alpha(x)=Q(x)[[r,r]]+\sum_{j}[P(x\cdot a_j)-P(x)P(a_j)](r_{12}-r_{21})\otimes b_j,\forall x\in A,\eqno (5.4)$$
where
$$[[r,r]]=r_{13}\cdot r_{12}-r_{23}\cdot r_{21}+[r_{23},r_{12}]-[r_{13},r_{21}]-[r_{13},r_{23}],\eqno (5.5)$$
and $Q(x)=L_x\otimes 1\otimes 1+1\otimes L_x\otimes 1+1\otimes 1\otimes {\rm ad}x$,
$P(x)=L_x\otimes 1+1\otimes L_x$ for any $x\in A$.

Before proving this result, let us explain the notations. Let
$(A,\cdot)$ be a left-symmetric algebra and $r=\sum_i a_i\otimes
b_i$,
set
$$r_{12}=\sum_ia_i\otimes b_i\otimes 1,r_{21}=\sum_i b_i\otimes a_i\otimes 1;\;
r_{13}=\sum_{i}a_i\otimes 1\otimes b_i;\;r_{23}=\sum_i1\otimes a_i\otimes b_i\;\in
U({\cal G}(A)),\eqno (5.6)$$
where $U({\cal G}(A))$ is the universal enveloping algebra of the sub-adjacent Lie algebra ${\cal G}(A)$. Set
\begin{eqnarray*}
&&r_{13}\cdot r_{12}=\sum_{i,j}a_i\cdot a_j\otimes b_j\otimes b_i;\;\;
r_{23}\cdot r_{21}=\sum_{i,j} b_j\otimes a_i\cdot a_j\otimes b_i;\\
&&[r_{23},r_{12}]=r_{23}\cdot r_{12}-r_{12}\cdot r_{23}=\sum_{ij} a_j\otimes [a_i,b_j]\otimes b_i;\\
&&[r_{13},r_{21}]=r_{13}\cdot r_{21}-r_{21}\cdot r_{13}=\sum_{ij} [a_i,b_j]\otimes a_j\otimes b_i;\\
&&[r_{13},r_{23}]=r_{13}\cdot r_{23}-r_{23}\cdot r_{13}=\sum_{ij}a_i\otimes a_j\otimes [b_i,b_j].
\end{eqnarray*}

{\bf Proof}\quad Let $x\in A$. After rearranging the terms suitably, we divide $J_\alpha(x)$ into three parts:
$$J_\alpha(x)=(F1)+(F2)+(F3),$$
where
\begin{eqnarray*}
(F1)&=&\sum_{i,j}
\{(x\cdot a_j)\cdot a_i\otimes b_i-b_i\otimes (x\cdot a_j)\cdot a_i+b_i\cdot a_j\otimes x\cdot a_i
-x\cdot a_i\otimes b_i\cdot a_j
\\
&\mbox{}&\;\;\;+a_i\otimes [x\cdot a_j,b_i- a_i\otimes [x,b_i]\cdot a_j+
[x,b_i]\cdot a_j\otimes a_i-[x\cdot a_j,b_i]\otimes a_i
\}\otimes b_j;\\
(F2)&=&\sum_{i,j}\{a_j\cdot a_i\otimes b_i+a_i\otimes[a_j,b_i]-b_i\otimes a_j\cdot a_i
-[a_j,b_i]\otimes a_i\}\otimes [x,b_j]\\
(F3)&=&\sum_{i,j}\{-x\cdot a_i\otimes a_j\otimes [b_i,b_j]-a_i\otimes a_j\otimes[[x,b_i],b_j]-a_i\otimes x\cdot a_j\otimes [b_i,b_j]\\
&\mbox{}&\;\;\;+a_i\otimes a_j\otimes [[x,b_j],b_i]]\}.
\end{eqnarray*}
On the other hand,
\begin{eqnarray*}
&&Q(x)(r_{13}\cdot r_{12})=\sum_{i,j}\{ x\cdot(a_j\cdot a_i)\otimes b_i\otimes b_j+a_j\cdot a_i\otimes x\cdot b_i\otimes b_j+a_j\cdot
a_i\otimes b_i\otimes [x,b_j]\};\\
&&Q(x)(-r_{23}\cdot r_{21})=\sum_{i,j} -\{x\cdot b_i\otimes a_j\cdot a_i\otimes b_j+b_i\otimes x\cdot (a_j\cdot a_i)\otimes b_j+
b_i\otimes a_j\cdot a_i\otimes [x,b_j]\};\\
&&Q(x)([r_{23},r_{12}])=\sum_{i,j} \{x\cdot a_i\otimes [a_j,b_i]\otimes b_j+a_i\otimes x\cdot [a_j,b_i]\otimes b_j+a_i\otimes [a_j,b_i]
\otimes [x,b_j]\};\\
&&Q(x)(-[r_{13},r_{21}])=\sum_{i,j}-\{x\cdot [a_j, b_i]\otimes a_i\otimes b_j+[a_j, b_i]\otimes x\cdot a_i\otimes b_j+[a_j,b_i]\otimes a_i\otimes [x,b_j]\};\\
&&Q(x)(-[r_{13},r_{23}])=\sum_{i,j}-\{x\cdot a_i\otimes a_j\otimes [b_i,b_j]+a_i\otimes x\cdot a_j\otimes [b_i,b_j]+a_i\otimes a_j\otimes
[x,[b_i,b_j]]\};\\
&&\sum_j [P(x\cdot a_j)-P(x)P(a_j)](r_{12}-r_{21})\otimes b_j=\sum_{i,j}\{ (x\cdot a_j)\cdot a_i\otimes b_i +a_i\otimes
(x\cdot a_j)\cdot b_i\\
&&\hspace{3cm}-x\cdot a_i\otimes a_j\cdot b_i
-a_j\cdot a_i\otimes x\cdot b_i-x\cdot (a_j\cdot a_i)\otimes b_i-a_i\otimes x\cdot (a_j\cdot b_i)\\
&&\hspace {3cm} -(x\cdot a_j)\cdot b_i\otimes a_i -b_i\otimes
(x\cdot a_j)\cdot a_i+x\cdot b_i\otimes a_j\cdot a_i+a_j\cdot b_i\otimes x\cdot a_i\\
&&\hspace {3cm} +x\cdot (a_j\cdot b_i)\otimes a_i+b_i\otimes x\cdot (a_j\cdot a_i)\}\otimes b_j
\end{eqnarray*}
After rearranging the terms suitably, the sum of the terms whose third component is $b_j$ in
the right hand side of equation (5.4) is
$$ (F1')=(F1a)+(F1b)+(F1c)+(F1d)+(F1e)+(F1f)+(F1g),$$
where
\begin{eqnarray*}
(F1a)&=&\sum_{i,j}\{ x\cdot (a_j\cdot a_i)+ (x\cdot a_j)\cdot a_i - x\cdot (a_j\cdot a_i)\}\otimes b_i
\otimes b_j=\sum_{i,j}(x\cdot a_j)\cdot a_i\otimes b_i\otimes b_j;\\
(F1b)&=&\sum_{i,j} b_i\otimes \{- x\cdot (a_j\cdot a_i)-(x\cdot a_j)\cdot a_i+x\cdot (a_j\cdot a_i)\}\otimes b_j
=\sum_{i,j}-b_i\otimes (x\cdot a_j)\cdot a_i\otimes b_j;\\
(F1c)&=&\sum_{i,j}\{-[a_j, b_i]\otimes x\cdot a_i+a_j\cdot b_i\otimes x\cdot a_i\}\otimes b_j=\sum_{i,j}b_i\cdot a_j\otimes
x\cdot a_i\otimes b_j;\\
(F1d)&=&\sum_{i,j}\{ x\cdot a_i\otimes [a_j,b_i]-x\cdot a_i\otimes a_j\cdot b_i\}\otimes b_j=
\sum_{i,j}-x\cdot a_i\otimes b_i\cdot a_j\otimes b_j;\\
(F1e)&=&\sum_{ij}\{a_i\otimes x\cdot [a_j,b_i]+a_i\otimes (x\cdot a_j)\cdot b_i-a_i\otimes x\cdot (a_j\cdot b_i)\}\otimes b_j\\
\mbox{}&=& \sum_{i,j} a_i\otimes \{ -x\cdot (b_i\cdot a_j)+(x\cdot a_j)\cdot b_i\}\otimes b_j
=\sum_{i,j}a_i\otimes \{ [x\cdot a_j,b_i]-[x,b_i]\cdot a_j\}\otimes b_j;\\
(F1f)&=&\sum_{i,j}\{ -x\cdot [a_j, b_i]\otimes a_i-(x\cdot a_j)\cdot b_i\otimes a_i+x\cdot (a_j\cdot b_i)\otimes a_i\}\otimes b_j\\
\mbox{}&=&\sum_{i,j} \{ x\cdot (b_i\cdot a_j)-(x\cdot a_j)\cdot b_i\}\otimes a_i\otimes b_j=
\sum_{i,j}\{ [x,b_i]\cdot a_j-[x\cdot a_j,b_i]\}\otimes a_i\otimes b_j;\\
(F1g)&=&\sum_{i,j} \{ a_j\cdot a_i\otimes x\cdot b_i-x\cdot b_i\otimes a_j\cdot a_i
-a_j\cdot a_i\otimes x\cdot b_i+x\cdot b_i\otimes a_j\cdot a_i\}\otimes b_j=0.
\end{eqnarray*}
Therefore $(F1')=(F1)$. Obviously, the sum of the terms whose third
component is $[x,b_j]$ in the right hand side of equation (5.4) is
just $(F2)$. Moreover, the sum of the other terms in the right hand
side of equation (5.4) is $Q(x)(-[r_{13},r_{23}])$, which precisely
equals to $(F3)$ by Jacobi identity in the sub-adjacent Lie algebra
${\cal G}(A)$. Hence equation (5.4) holds.\hfill $\Box$

With the discussion above together, we have the following result.

{\bf Theorem 5.4}\quad Let $A$ be a left-symmetric algebra and $r\in
A\otimes A$. Then the map $\alpha$ defined by equation (5.1) induces
a left-symmetric algebra structure on $A^*$ such that $(A,A^*)$ is a
left-symmetric bialgebra if and only if the following two conditions
are satisfied:

(a) $[P(x\cdot y)-P(x)P(y)](r_{12}-r_{21})=0$ for any $x,y\in A$;

(b) $Q(x)[[r,r]]=0$,

\noindent where $[[r,r]]$ is given by equation (5.5) and $Q(x)=L_x\otimes 1\otimes 1+1\otimes L_x\otimes 1+1\otimes 1\otimes {\rm ad}x$,
$P(x)=L_x\otimes 1+1\otimes L_x$ for any $x\in A$.

A direct application of Theorem 5.4 is given as follows (cf. Example 4.3).

{\bf Theorem 5.5}\quad Let $(A,A^*,\alpha,\beta)$ be a
left-symmetric bialgebra. Then there is a canonical left-symmetric
bialgebra structure on $A\oplus A^*$ such that both the inclusions
$i_1:A\rightarrow A\oplus A^*$ and $i_2:A^*\rightarrow A\oplus A^*$
into the two summands are homomorphisms of left-symmetric
bialgebras.

{\bf Proof}\quad Let $r\in A\otimes A^*\subset (A\oplus A^*)\otimes
(A\oplus A^*)$ correspond to the identity map $id:A\rightarrow A$.
Let $\{ e_1,\cdots,e_n\}$ be a basis of $A$ and $\{ e_1^*,\cdots,
e_n^*\}$ be its dual basis. Then $r=\sum_{i} e_i\otimes e_i^*$.
Suppose that the left-symmetric algebra structure ``$*$'' on
$A\oplus A^*$ is given by ${\cal SD}(A)=A\bowtie^{{\rm
ad}^*_\cdot,-R_\cdot^*}_{{\rm ad}^*_\circ,-R^*_\circ}A^*$. Then by
Theorem 3.5, we have
$$x*y=x\cdot y,\;\;a^**b^*=a^*\circ b^*,\;x*a^*={\rm ad}_\cdot^*(x)a^*-R_\circ^*(a^*)x,\;
a^**x={\rm ad}_\circ (a^*)x-R^*_\cdot(x)a^*,$$
for any $x,y\in A, a^*,b^*\in A^*$.
Next we prove that $r$ satisfies the two conditions in Theorem 5.4. If so, then
$$\alpha_{{\cal SD}}(u)=(L_u\otimes 1+1\otimes {\rm ad}u)(r),\;\;\forall u\in {\cal SD}(A)$$
can induce a left-symmetric bialgebra structure on ${\cal SD}(A)$.

For any $\lambda,\mu\in {\cal SD}(A)$, equation (5.2) is equivalent to
\begin{eqnarray*}
&&\sum_i\{(\lambda*\mu)*e_i\otimes e_i^*-\lambda*(\mu*e_i)\otimes e_i^*-(\lambda*\mu)*e_i^*\otimes e_i+\lambda*(\mu*e_i^*)\otimes e_i\\
&&+e_i\otimes (\lambda*\mu)*e_i^*-e_i\otimes \lambda*(\mu*e_i^*)+e_i^*\otimes \lambda*(\mu*e_i)-e_i^*\otimes (\lambda*\mu)*e_i\\
&&-\lambda*e_i\otimes \mu*e_i^*+\mu*e_i^*\otimes \lambda*e_i+\lambda*e_i^*\otimes \mu*e_i-\mu*e_i\otimes \lambda*e_i^*\}=0
\end{eqnarray*}
We can prove the equation above in the following four cases: (I)
$\lambda,\mu\in A$; (II) $\lambda,\mu\in A^*$; (III) $\lambda\in
A,\mu\in A^*$ and (IV) $\lambda\in A^*,\mu\in A$. As an example, we
give a detailed proof of the first case (the proof of other cases is
similar). Let $\lambda=e_k$, $\mu=e_l$, then we have
\begin{eqnarray*}
\sum_i e_i\otimes (e_k*e_l)*e_i^*&=&\sum_i -[e_k\cdot e_l,e_i]\otimes e_i^*+\sum_{i,m}\langle e_m^*\circ e_i^*, e_k\cdot e_l\rangle   e_i\otimes e_m;\\
\sum_i -(e_k*e_l)*e_i^*\otimes e_i&=&\sum_i e_i^*\otimes [e_k\cdot e_l, e_i]-\sum_{i,m} \langle e_m^*\circ e_i^*, e_k\cdot e_l\rangle  e_m\otimes e_i;\\
\sum_i -e_i\otimes e_k*(e_l*e_i^*)&=&\sum_i -[e_l,[e_k,e_i]]\otimes e_i^*-
\sum_{i,m} \langle e_m^*\circ e_i^*, e_l\rangle   e_i\otimes e_k\cdot e_m\\
&+&\sum_{i,m,n}\langle [e_l,e_m],e_i^*\rangle  \langle e_n^*\circ e_m^*, e_k\rangle   e_i\otimes e_n;\\
\sum_i e_k*(e_l*e_i^*)\otimes e_i&=&\sum_i e_i^*\otimes [e_l,[e_k,e_i]]+
\sum_{i,m} \langle e_m^*\circ e_i^*, e_l\rangle  e_k\cdot e_m\otimes e_i\\
&-&\sum_{i,m,n}\langle [e_l,e_m],e_i^*\rangle  \langle e_n^*\circ e_m^*,e_k\rangle  e_n\otimes e_i;\\
\sum_i-e_k*e_i\otimes e_l* e_i^*&=& \sum_ie_k\cdot [e_l,e_i]\otimes e_i^*-\sum_{i,m}\langle e_m^*\circ e_i^*,e_l\rangle  e_k\cdot e_i\otimes e_m;\\
\sum_i e_l*e_i^*\otimes e_k*e_i&=&\sum_i -e_i^* \otimes  e_k\cdot [e_l,e_i]+\sum_{i,m}\langle e_m^*\circ e_i^*, e_l\rangle  e_m\otimes e_k\cdot e_i;\\
\sum_i-e_l*e_i\otimes e_k* e_i^*&=& \sum_ie_l\cdot [e_k,e_i]\otimes e_i^*-\sum_{i,m}\langle e_m^*\circ e_i^*,e_k\rangle  e_l\cdot e_i\otimes e_m;\\
\sum_i e_k*e_i^*\otimes e_l*e_i&=&\sum_i -e_i^* \otimes  e_l\cdot
[e_k,e_i]+\sum_{i,m}\langle e_m^*\circ e_i^*, e_k\rangle  e_m\otimes
e_l\cdot e_i.
 \end{eqnarray*}
The sum of the terms which are in $A\otimes A^*$ or $A^*\otimes A$
is zero since
\begin{eqnarray*}
&& (e_k\cdot e_l)\cdot e_i -e_k\cdot (e_l\cdot e_i)-[e_k\cdot e_l, e_i]-[e_l,[e_k,e_i]]
+e_k\cdot [e_l,e_i]+e_l\cdot [e_k,e_i]\\
&& =e_i\cdot (e_k\cdot e_l)-(e_i\cdot e_k)\cdot e_l-e_k(e_i\cdot e_l)+(e_k\cdot e_i)\cdot e_l=0.
\end{eqnarray*}
The coefficient of $e_i\otimes e_m$ in the sum of the terms which
are in $A\otimes A$ is
\begin{eqnarray*}
&&\langle e_m^*\circ e_i^*, e_k\cdot e_l\rangle  -\langle e_i^*\circ e_m^*,e_k\cdot e_l\rangle  \\
&&\hspace{1cm} -\sum_n\{ \langle e_n^*\circ e_i^*,e_l\rangle  \langle e_m^*, e_k\cdot e_n\rangle  +\langle [e_l,e_n],e_i^*\rangle  \langle e_m^*\circ e_n^*,e_k\rangle  \\
&&\hspace{1cm} \langle e_i^*,e_k\cdot e_n\rangle  \langle e_n^*\circ e_m^*,e_l\rangle  -\langle e_m^*,[e_l,e_n]\rangle  \langle e_i^*\circ e_n^*,e_k\rangle  \\
&&\hspace{1cm}-\langle e_i^*,e_k\cdot e_n\rangle  \langle e_m^*\circ e_n^*,e_l\rangle  +\langle e_m^*,e_k\cdot e_n\rangle  \langle e_i^*\circ e_n^*,e_l\rangle  \\
&&\hspace{1cm}-\langle e_i^*,e_l\cdot e_n\rangle  \langle e_m^*\circ e_n^*,e_k\rangle  +\langle e_m^*,e_l\cdot e_n\rangle  \langle e_i^*\circ e_n^*,e_k\rangle  \}\\
&&=\langle -L^*_\cdot (e_k)[e_m^*,e_i^*],e_l\rangle  +\langle
L_\cdot^*(e_k)e_m^*\circ e_i^*,e_l\rangle
+\langle -L_\cdot^*(L_\circ^*(e_m^*)e_k)e_i^*, e_l\rangle  \\
&&\hspace{1cm} -\langle L_\cdot^*(e_k)e_i^*\circ e_m^*,e_l\rangle  +\langle L_\cdot^*(L_\circ^*(e_i^*)e_k)e_m^*, e_l\rangle  \\
&&\hspace{1cm}+\langle e_m^*\circ L_\cdot^*(e_k)e_i^*,e_l\rangle
-\langle e_i^*\circ L_\cdot^*(e_k)e_m^*,e_l\rangle  =0
\end{eqnarray*}
since $({\cal G}(A),{\cal G}(A^*),L_\cdot^*,L_\circ^*)$ is a matched pair of Lie algebras. Hence
$P(e_k\cdot e_l)-P(e_k)P(e_l)(r_{12}-r_{21})=0$. Furthermore,
\begin{eqnarray*}
[[r,r]]&=& \sum_{i,j}\{ e_i\cdot e_j\otimes e_j^*\otimes e_i^*+e_j\otimes [e_i,e_j^*]\otimes e_i^*
-e_j^*\otimes e_i\cdot e_j\otimes e_i^*\\
&\mbox{}&-[e_i,e_j^*]\otimes e_j\otimes e_i^*-e_i\otimes e_j\otimes [e_i^*,e_j^*]\}
\end{eqnarray*}
Note that
\begin{eqnarray*}
&&\sum_{i,j}\{ e_i\cdot e_j\otimes e_j^*\otimes e_i^*+e_j\otimes L_\cdot^* (e_i)e_j^*\otimes e_i^*\}=0;\\
&&\sum_{i,j}\{-e_j\otimes L_\circ^*(e_j^*)e_i\otimes e_i^*+L_\circ^*(e_j^*)e_i\otimes e_j\otimes e_i^*
-e_i\otimes e_j\otimes [e_i^*,e_j^*]\}=0;\\
&&\sum_{i,j}\{-L_\cdot^*(e_i)e_j^*\otimes e_j\otimes e_i^*-e_j^*\otimes e_i\cdot e_j\otimes e_i^*\}=0.
\end{eqnarray*}
Therefore $[[r,r]]=0$. Hence ${\cal SD}(A)$ is a left-symmetric bialgebra.

For $e_i\in A$, we have
\begin{eqnarray*}
\alpha_{SD}(e_i)&=& \sum_{j}\{ e_i\cdot e_j\otimes e_j^*+e_j\otimes [e_i,e_j^*]\}\\
&=&\sum_{j,m}\{ e_i\cdot e_j\otimes e_j^*-\langle e_j^*,e_i\cdot e_m\rangle  e_j\otimes e_m^*+\langle e_j^*\circ e_m^*,e_i\rangle  e_j\otimes e_m\}\\
&=&\sum_{j,m}\langle e_j^*\circ e_m^*,e_i\rangle  e_j\otimes e_m\\
&=& \alpha_A (e_i)
\end{eqnarray*}
Therefore the inclusion $i_1:A\rightarrow A\oplus A^*$ is a homomorphism of left-symmetric bialgebras. Similarly,
the inclusion $i_2:A^*\rightarrow A\oplus A^*$ is also a homomorphism of left-symmetric bialgebras since
$\alpha_{SD}(e_i^*)=\alpha_{A^*}(e_i^*)$.\hfill $\Box$

{\bf Definition 5.2}\quad Let $(A,A^*)$ be a left-symmetric bialgebra. With the left-symmetric bialgebra structure
given in Theorem 5.5, $A\oplus A^*$ is called a symplectic double of $A$. We denote it by ${\cal SD}(A)$.

Therefore, by Theorem 5.5, we have the following conclusion.

{\bf Proposition 5.6}\quad  Let $(A,A^*)$ be a left-symmetric
bialgebra. Then the symplectic double ${\cal SD}(A)$ of $A$ is a
left-symmetric bialgebra and its sub-adjacent Lie algebra is a
parak\"ahler Lie algebra with the symplectic form $\omega_p$ given
by equation (2.16).

At the end of this section, we would like to point out that, unlike the symmetry of 1-cocycle of ${\cal G}(A)$ and ${\cal G}(A^*)$ appearing in the definition of a left-symmetric bialgebra $(A,A^*)$, it is not necessary that $\beta$ is also a 1-coboundary of ${\cal G}(A^*)$ for a coboundary left-symmetric bialgebra $(A,A^*,\alpha,\beta)$, where $\alpha$ is given by equation (5.1).

\section{$S$-equation and its properties}

The simplest way to satisfy the two conditions of Theorem 5.4 is to assume that $r$ is symmetric and $[[r,r]]=0$.
Note that
$$[[r,r]]=r_{21}\cdot r_{13}-r_{12}\cdot r_{23}-[r_{13},r_{23}]+r_{13}\cdot(r_{12}-r_{21})+r_{23}
\cdot (r_{12}-r_{21}).\eqno (6.1)$$

Therefore, by Theorem 5.4, we have the following conclusion.

{\bf Proposition 6.1}\quad Let $A$ be a left-symmetric algebra and
$r\in A\otimes A$. Suppose $r$ is symmetric. Then the map $\alpha$
defined by equation (5.1) induces a left-symmetric algebra structure
on $A^*$ such that $(A,A^*)$ is a left-symmetric bialgebra if
$$[[r,r]]=-r_{12}\cdot r_{13}+r_{12}\cdot r_{23}+[r_{13},r_{23}]=0.\eqno (6.2)$$

{\bf Definition 6.1}\quad Let $A$ be a left-symmetric algebra and
$r\in A\otimes A$. Then equation (6.2) is called ${S}$-equation in
$A$.

{\bf Remark 1}\quad Let the products on $\{1,2,3\}$ correspond to
the products of $A$: for example, for any $i,j,k=1,2,3$, $(i\cdot
j)\cdot k$ or $i\cdot (j\cdot k)$ corresponds to the (order) product
$(x\cdot y)\cdot z$ or $x\cdot (y\cdot z)$ of $A$ respectively, and
so on. For any symmetric $r\in A\otimes A$, set (for any $r_{ij},
i<j, i,j=1,2,3$)
$$r_{ij}\cdot r_{jl}\longleftrightarrow j\cdot (i\cdot l),\;\;
r_{ij}\cdot r_{il}\longleftrightarrow i\cdot (j\cdot l),\;\;
r_{ij}\cdot r_{lj}\longleftrightarrow (i\cdot l)\cdot j,\;\;
r_{ij}\cdot r_{li}\longleftrightarrow (j\cdot l)\cdot i.$$ Then the
$S$-equation (6.2) in a left-symmetric algebra  corresponds to the
``left-symmetry'' of the products. It is similar to the relation
between the classical Yang-Baxter equation in a Lie algebra ${\cal
G}$ and the Jacobi identity of ${\cal G}$ ([BD],[D],[Se]).

{\bf Remark 2}\quad  Let $\sigma_{123},\sigma_{132}:A\otimes
A\otimes A\rightarrow A\otimes A\otimes A$ be two linear maps
satisfying $$\sigma_{123}(x\otimes y\otimes z)=z\otimes x\otimes
y,\;\; \sigma_{132}(x\otimes y\otimes z)=y\otimes z\otimes
x,\;\;\forall \; x,y,z\in A.$$ Suppose that $r$ is symmetric. Then
it is easy to know that under the action of $\sigma_{123}$ and
$\sigma_{132}$ respectively, the $S$-equation (6.2) turns to be
$$[r_{12},r_{13}]-r_{23}\cdot (r_{12}-r_{13})=0;$$
$$[r_{12},r_{23}]-r_{13}\cdot (r_{12}-r_{23})=0,$$
respectively.

Let $A$ be a vector space. For any $r\in A\otimes A$, $r$ can be regarded as a map
from $A^*$ to $A$ in the following way:
$$\langle u^*\otimes v^*,r\rangle  =\langle u^*,r(v^*)\rangle  ,\;\;\forall u^*,v^*\in A^*.\eqno (6.3)$$

{\bf Proposition 6.2}\quad Let $(A,\cdot)$ be a left-symmetric and
$r\in A\otimes A$ be a symmetric solution of $S$-equation in $A$.
Then the left-symmetric algebra and its sub-adjacent Lie algebra
structure on the symplectic double ${\cal SD}(A)$ can be given from
the products in $A$ as follows:

(a) $a^**b^*=a^*\circ b^*=-R^*_\cdot (r(b^*))a^*+{\rm ad}_\cdot^*(r(a^*))b^*$, for any $a^*,b^*\in A^*$;\hfill (6.4)

(b) $[a^*,b^*]=a^*\circ b^*-b^*\circ a^*=L^*_\cdot (r(a^*))b^*-L_\cdot^*(r(b^*))a^*$, for any $a^*,b^*\in A^*$;\hfill (6.5)

(c) $x*a^*=x\cdot r(a^*)-r({\rm ad}_\cdot^*(x)a^*) + {\rm ad}_\cdot^* (x)a^*$, for any $x\in A$, $a^*\in A^*$;\hfill (6.6)

(d) $a^**x=r(a^*)\cdot x+r(R_\cdot^* (x) a^*)-R_\cdot^*(x)a^*$, for any $x\in A$, $a^*\in A^*$;\hfill (6.7)

(e) $[x,a^*]=[x,r(a^*)]-r(L_\cdot^*(x)a^*)+L_\cdot^*(x)a^*$, for any $x\in A$, $a^*\in A^*$.\hfill (6.8)

{\bf Proof}\quad Let $\{ e_1,\cdots,e_n\}$ be a basis of $A$ and $\{ e_1^*,\cdots, e_n^*\}$ be its dual basis.
Suppose that $e_i\cdot e_j=\sum_{i,j}c_{ij}^ke_k$ and
$r=\sum_{i,j}a_{ij}e_i\otimes e_j$, where $a_{ij}=a_{ji}$. Then from the proof of Proposition 5.1, we know that (for any $k,l$)
\begin{eqnarray*}
e_k^*\circ e_l^*&=&\sum_{i,t}\{a_{tl}c_{it}^k+a_{kt}(c_{it}^l-c_{ti}^l)\}e_i^*=
\sum_{i,t}\{ a_{tl}\langle e_i\cdot e_t,e_k^*\rangle  +a_{kt}\langle [e_i,e_t],e_l^*\rangle  \}e_i^*\\
&=& \sum_{i}\{ \langle e_i\cdot r(e_l^*),e_k^*\rangle  +\langle [e_i,r(e_k^*)],e_l^*\rangle  \}e_i^*\\
&=&-R^*_\cdot (r(e_l^*))e_k^*+{\rm ad}_\cdot^*(r(e_k^*))e_l^*.
\end{eqnarray*}
Hence equation (6.4) holds. Therefore, we have
\begin{eqnarray*}
-R_\circ^*(e_k^*)e_l&=& \sum_i \langle -R_\circ^*(e_k^*)e_l,e_i^*\rangle  e_i=\sum_i \langle e_l,e_i^*\circ e_k^*\rangle  e_i\\
&=&\sum_i\langle e_l, -R^*_\cdot (r(e_k^*))e_i^*+{\rm ad}_\cdot^*(r(e_i^*))e_k^*\rangle   e_i\\
&=&e_l\cdot r(e_k^*)-\sum_i\langle [r(e_i^*),e_l],e_k^*\rangle  e_i\\
&=&e_l\cdot r(e_k^*)-r({\rm ad}_\cdot^* (e_l) e_k^*).
\end{eqnarray*}
Since $e_l*e_k^*={\rm ad}_\cdot^*(e_l)e_k^*-R_\circ^*(e_k^*)e_l$, we get equation (6.6). Similarly, we can prove equation (6.7). Equations (6.5) and (6.8) then follow immediately.\hfill $\Box$

{\bf Definition 6.2}\quad Let $(A,\cdot)$ be a left-symmetric
algebra. A bilinear form ${\cal B}:A\otimes A\rightarrow {\bf F}$ is
called a 2-cocycle of $A$ if
$${\cal B}(x\cdot y,z)-{\cal B}(x,y\cdot z)={\cal B}(y\cdot x,z)-{\cal B}(y,x\cdot z),\forall x,y,z\in A.\eqno (6.9)$$

In fact, the above definition is precisely the definition of a
2-cocycle of a left-symmetric algebra into the trivial bimodule
${\bf F}$ ([SW]). It is equivalent to the following central
extension: there exists a left-symmetric algebra structure on
$A\oplus {\bf F}c$ given by
$$x*y=x\cdot y+{\cal B}(x,y)c,\;x*c=c*x=c*c=0,\;\;\forall x,y\in A, \eqno (6.10)$$
if and only if ${\cal B}$ is a 2-cocycle of $A$. Furthermore, for
any 2-cocycle ${\cal B}$ of a left-symmetric algebra $A$, it is easy
to know that ([Ku2])
$$\omega (x,y)={\cal B}(x,y)-{\cal B}(y,x),\;\;\forall x,y\in A,\eqno (6.11)$$
is a 2-cocycle of the sub-adjacent Lie algebra ${\cal G}(A)$.

{\bf Theorem 6.3}\quad Let $A$ be a left-symmetric algebra and $r\in
A\otimes A$. Suppose that $r$ is symmetric and nondegenerate. Then
$r$ is a solution of $S$-equation in $A$ if and only if the inverse
of the isomorphism $A^*\rightarrow A$ induced by $r$, regarded as a
bilinear form ${\cal B}$ on $A$, is a 2-cocycle of $A$. That is,
${\cal B}(x,y)=\langle r^{-1}x,y\rangle  $ for any $x,y\in A$.

{\bf Proof}\quad Let $r=\sum_i a_i\otimes b_i$. Since $r$ is
symmetric, we have $\sum_i a_i\otimes b_i=\sum_ib_i\otimes a_i$.
Therefore $r(v^*)=\sum_i\langle v^*,a_i\rangle  b_i=\sum_i\langle
v^*,b_i\rangle  a_i$ for any $v^*\in A^*$. Since $r$ is
nondegenerate, for any $x,y,z\in A$, there exist $u^*,v^*,w^*\in
A^*$ such that $x=r(u^*),y=r(v^*),z=r(w^*)$. Therefore
\begin{eqnarray*}
{\cal B}(x\cdot y,z)&=& \langle r(u^*)\cdot r(v^*), w^*\rangle  =\sum_{i,j}\langle u^*,b_i\rangle  \langle v^*,b_j\rangle  \langle w^*,a_i\cdot a_j\rangle  \\
&=&\langle u^*\otimes v^*\otimes w^*, r_{13}\cdot r_{23}\rangle  ;\\
-{\cal B}(x,y\cdot z) &=&-\langle u^*, r(v^*)\cdot r(w^*)\rangle  =\sum_{i,j}-\langle v^*,b_i\rangle  \langle w^*,b_j\rangle  \langle u^*,a_i\cdot a_j\rangle  \\
&=&\langle u^*\otimes v^*\otimes w^*, -r_{12}\cdot r_{13}\rangle  ;\\
-{\cal B}(y\cdot x,z)&=& -\langle r(v^*)\cdot r(u^*), w^*\rangle  =-\sum_{i,j}\langle v^*,b_i\rangle  \langle u^*,b_j\rangle  \langle w^*,a_i\cdot a_j\rangle  \\
&=&\langle u^*\otimes v^*\otimes w^*, -r_{23}\cdot r_{13}\rangle  ;\\
{\cal B}(y,x\cdot z) &=&\langle v^*, r(u^*)\cdot r(w^*)\rangle  =\sum_{i,j}\langle u^*,b_i\rangle  \langle w^*,b_j\rangle  \langle v^*,a_i\cdot a_j\rangle  \\
&=&\langle u^*\otimes v^*\otimes w^*, r_{12}\cdot r_{23}\rangle  .
\end{eqnarray*}
Therefore ${\cal B}$ is a (symmetric) 2-cocycle of $A$ if and only
if $\langle u^*\otimes v^*\otimes w^*,[[r,r]]\rangle  =0$ for any
$u^*,v^*,w^*\in A^*$, if and only if $[[r,r]]$=0.\hfill $\Box$

{\bf Corollary 6.4}\quad Let $(A,\cdot)$ be a left-symmetric and
$r\in A\otimes A$ be a nondegenerate symmetric solution of
$S$-equation in $A$. Suppose the left-symmetric algebra structure
$``\circ"$ on $A^*$ is induced by $r$ from equation (6.4). Then we
have
$$a^*\circ b^*=r^{-1}(r(a^*)\cdot r(b^*)),\;\;\forall a^*,b^*\in A^*.\eqno (6.12)$$
Therefore $r:A^*\rightarrow A$ is an isomorphism of left-symmetric algebras.

{\bf Proof}\quad For any $x,y\in A$, set ${\cal B}(x,y)=\langle
r^{-1}(x),y\rangle  $. Then ${\cal B}$ is a 2-cocycle of $A$. For
any $a^*,b^*\in A^*$ and $x\in A$, from Proposition 6.2 and the
above theorem, we know that
\begin{eqnarray*}
\langle a^*\circ b^*,x\rangle  &=&\langle x\cdot r(b^*), a^*\rangle  -\langle [r(a^*),x],b^*\rangle  \\
&=&{\cal B}(x\cdot r(b^*),r(a^*))-{\cal B}([r(a^*),x],r(b^*))={\cal B}(x,r(a^*)\cdot r(b^*))\\
&=& \langle r^{-1}(r(a^*)\cdot r(b^*)),x\rangle  .
\end{eqnarray*}
Hence equation (6.12) holds and therefore $r$ is an isomorphism of
left-symmetric algebras. \hfill $\Box$

{\bf Definition 6.3}\quad A left-symmetric algebra $A$ over the real field ${\bf R}$ is called Hessian if there exists
a symmetric and positive definite 2-cocycle ${\cal B}$ of $A$, that is, there exists a 2-cocycle ${\cal B}$ of $A$ which is an inner product.

In geometry, a Hessian manifold $M$ is a flat affine manifold provided with a Hessian metric $g$, that is, $g$ is a Riemannian metric such that for any each point $p\in M$ there exists a $C^\infty$-function $\varphi$ defined on a neighborhood of
$p$ such that $g_{ij}=\frac{\partial^2\varphi}{\partial x^i\partial x^j}$. A Hessian left-symmetric algebra corresponds to an affine Lie group $G$ with a $G$-invariant Hessian metric ([Sh]).

Hence, we have the following conclusion.

{\bf Proposition 6.5}\quad Let $(A,{\cal B})$ be a Hessian
left-symmetric algebra with the inner product ${\cal B}$. Then there
exists a basis $\{ e_1,\cdots,e_n\}$ such that ${\cal
B}(e_i,e_j)=\delta_{ij}$. Suppose $e_i\cdot e_j=\sum_kc_{ij}^ke_k$.
Then
$$ c_{ij}^k-c_{ji}^k+c_{ik}^j-c_{jk}^i=0,\;\;\forall i,j,k.\eqno (6.13)$$
Under this basis, the corresponding symmetric solution of $S$-equation in $A$ is given by
$$ r=\sum_{i}^n e_i\otimes e_i.\eqno (6.14)$$
Moreover, let $\{e_1^*,\cdots, e_n^*\}$ be the dual basis of
$\{e_1,\cdots,e_n\}$. Then the left-symmetric algebra and its
sub-adjacent Lie algebra structures on the symplectic double ${\cal
SD}(A)$ are given by
$$e_i^**e_j^*=r(e_i\cdot e_j)=\sum_k^n c_{ij}^ke_k^*;\;[e_i^*,e_j^*]=r([e_i,e_j])=\sum_k^n(c_{ij}^k-c_{ji}^k)e_k^*;$$
$$e_i*e_j^*=\sum_k^n [c_{kj}^ie_k+(c_{ki}^j-c_{ik}^j)e_k^*];\;
e_j^**e_i=\sum_k^n[(c_{ji}^k-c_{ij}^k)e_k+c_{ki}^je_k^*];\;
[e_i,e_j^*]=\sum_k^n(c_{jk}^ie_k-c_{ik}^je_k^*).\eqno (6.15)$$ On
the other hand, if equation (6.14) is a symmetric solution of
$S$-equation in a left-symmetric algebra $A$, then the structure
constants $\{c_{ij}^k\}$ associated to the basis $\{
e_1,\cdots,e_n\}$ satisfy equation (6.13).

{\bf Example 6.1}\quad For an arbitrary nonassociative algebra $(A,\cdot)$, there is an
 ``invariant" bilinear form ${\cal B}$ defined as a trace form ([Sc]):
$${\cal B}(x\cdot y,z)={\cal B}(x,y\cdot z),\;\;{\rm or}\;\;{\rm equivalently},
{\cal B}(L_x(y),z)={\cal B}(x,R_y(z)),\;\forall x,y,z\in A.\eqno
(6.16)$$ Obviously, for a left-symmetric algebra $A$, a trace form
is a 2-cocycle of $A$. Hence a nondegenerate symmetric trace form of
a left-symmetric algebra $A$ can give a symmetric solution of
$S$-equation in $A$. For example, ${\bf F}$ can be regarded as a
1-dimensional associative algebra with a basis $\{e\}$ satisfying
$e\cdot e=e$. Then $r=e\otimes e$ as a solution of $S$-equation in
${\bf F}$ corresponds to the trace form ${\cal B}(e,e)=1$. Moreover,
a left-symmetric algebra with a nondegenerate symmetric trace form
must be associative and hence it is a Frobenius algebra ([Bo2]).

{\bf Example 6.2}\quad Let $(A,\cdot)$ be a left-symmetric algebra. There is another ``invariant" bilinear
form ${\cal B}$ in the following sense ([BM4]):
$${\cal B}(x\cdot y,z)=-{\cal B}(y,x\cdot z),\;\;{\rm or}\;\;{\rm equivalently},
 {\cal B}(L_x(y),z)=-{\cal B}(y,L_x(z)),\;\forall x,y,z\in A.\eqno (6.17)$$
The existence of such a nondegenerate bilinear form on $A$ is
equivalent to the fact that $L$ is isomorphic to $L^*$ as
representations of ${\cal G}(A)$. Obviously, such a nondegenerate
symmetric bilinear form on a left-symmetric algebra $A$ is also a
2-cocycle of $A$. In fact, we have constructed some examples of
parak\"ahler Lie algebras from these bilinear forms in [Bai2]. Here,
we explain why such a bilinear form appears in our study. As it was
pointed in [BM4], such bilinear forms are degenerate on most of
left-symmetric algebras.

Next we turn to the general symmetric solutions of $S$-equation.

{\bf Theorem 6.6}\quad Let $(A,\cdot)$ be a left-symmetric algebra
and $r\in A\otimes A$ be symmetric. Then $r$ is a solution of
$S$-equation in $A$ if and only if $r$ satisfies
$$[r(a^*),r(b^*)]=r(L_\cdot^*(r(a^*))b^*-L_\cdot^*(r(b^*))a^*),\;\;\forall a^*,b^*\in A^*.\eqno (6.18)$$

{\bf Proof}\quad Let $\{ e_1,\cdots,e_n\}$ be a basis of $A$ and $\{ e_1^*,\cdots,e_n^*\}$ be its dual basis. Suppose that $e_i\cdot e_j=\sum_k c_{ij}^k e_k$ and $r=\sum_{i,j}a_{ij}e_i\otimes e_j$, $a_{ij}=a_{ji}$. Hence
$r(e_i^*)=\sum_k a_{ik}e_k$. Then $r$ is a solution of $S$-equation in $A$ if and only if (for any $i,j,k$)
$$\sum_{t,l}\{ -c_{tl}^ia_{tj}a_{lk}+c_{tl}^ja_{it}a_{lk}+(c_{tl}^k-c_{lt}^k)a_{it}a_{lj}\}=0.$$
The left-hand side of the above equation is precisely the
coefficient of $e_k$ in
$$[r(e_i^*),r(e_j^*)]-r(L_\cdot^*(r(e_i^*))e_j^*-L_\cdot^*(r(e_j^*))e_i^*).$$
Therefore the conclusion follows.\hfill $\Box$

{\bf Definition 6.4}\quad ([Ku3])\quad Let ${\cal G}$ be a Lie
algebra and $\rho:{\cal G}\rightarrow gl(V)$ be its representation.
A linear map $T:V\rightarrow {\cal G}$ is called an ${\cal
O}$-operator associated to $\rho$ if $T$ satisfies
$$[T(u), T(v)]=T(\rho(T(u))v-\rho(T(v))u),\forall u,v\in V.\eqno (6.19)$$

Obviously, for a left-symmetric algebra $A$, the identity map $id$ is an ${\cal O}$-operator associated to the regular representation of the sub-adjacent Lie algebra ${\cal G}(A)$. On the other hand,

{\bf Lemma 6.7}\quad ([Bai3])\quad Let ${\cal G}$ be a Lie algebra
and $\rho:{\cal G}\rightarrow gl(V)$ be its representation.  Let
$T:V\rightarrow {\cal G}$ be an ${\cal O}$-operator associated to
$\rho$. Then the product
$$u\circ v=\rho(T(u))v,\;\;\forall u,v\in V\eqno (6.20)$$
defines a left-symmetric algebra structure on $V$.  Therefore $V$ is
a Lie algebra as the sub-adjacent Lie algebra of this left-symmetric
algebra and $T$ is a homomorphism of Lie algebras. Furthermore,
$T(V)=\{T(v)|v\in V\}\subset {\cal G}$ is a Lie subalgebra of ${\cal
G}$ and there is an induced left-symmetric algebra structure on
$T(V)$ given by
$$T(u)\cdot T(v)=T(u\circ v)=T(\rho (T(u))v),\;\;\forall u,v\in V.\eqno (6.21)$$
Moreover, its sub-adjacent Lie algebra structure is
just the Lie subalgebra structure of ${\cal G}$
and $T$ is a homomorphism of left-symmetric algebras.

Therefore we have the following result.

{\bf Corollary 6.8}\quad Let $(A,\cdot)$ be a left-symmetric
algebra. Let $r\in A\otimes A$ be a symmetric solution of
$S$-equation in $A$. Then $r$ is an ${\cal O}$-operator associated
to $L_\cdot^*$. Therefore there is a left-symmetric algebra
structure on $A^*$ given by
$$a^*\circ 'b^*=L_\cdot^*(r(a^*))b^*,\;\;\forall a^*,b^*\in A^*.\eqno (6.22)$$
It has the same sub-adjacent Lie algebra of the left-symmetric
algebra in $A^*$ given by equation (6.4), which is induced by $r$ in
the sense of coboundary left-symmetric bialgebras. If $r$ is
nondegenerate, then there is a new compatible left-symmetric algebra
structure on ${\cal G}(A)$ given by
$$x\cdot'y= r(L_\cdot^*(x)r^{-1}y),\;\;\forall x,y\in A,\eqno (6.23)$$
which is just the left-symmetric algebra structure given by
$${\cal B}(x\cdot'y,z)=-{\cal B}(y,x\cdot z),\forall x,y,z\in A,\eqno (6.24)$$
where ${\cal B}$ is the symmetric nondegenerate 2-cocycle of $A$
induced by $r^{-1}$, that is, for any $x,y\in A$, ${\cal
B}(x,y)=\langle r^{-1}(x),y\rangle  $.

{\bf Remark}\quad With the above conditions and the left-symmetric
algebra structure on $A^*$ given by equation (6.22), $(A,A^*)$ is a
left-symmetric bialgebra if and only if the following two conditions
hold: (for any $x,y\in A$ and $a^*,b^*\in A^*$)
$$ L_\cdot^*[x\cdot r(a^*)-r(L_\cdot^*(x)a^*)]b^*=
L_\cdot^*[x\cdot r(b^*)-r(L_\cdot^*(x)b^*)]a^*;\eqno (6.25)$$
$$[x\cdot r(a^*)-r(L_\cdot^*(x)a^*)]\cdot y=
[y\cdot r(a^*)-r(L_\cdot^*(y)a^*)]\cdot x.\eqno (6.26)$$

{\bf Theorem 6.9}\quad Let ${\cal G}$ be a Lie algebra. Let
$\rho:{\cal G}\rightarrow gl(V)$ be a representation of ${\cal G}$
and $\rho^*:{\cal G}\rightarrow gl(V^*)$ be its dual representation.
Suppose that $T:V\rightarrow {\cal G}$ is an ${\cal O}$-operator
associated to $\rho$. Then
$$r=T+T^{21}\eqno (6.27)$$
is a symmetric solution of $S$-equation in
$T(V)\ltimes_{\rho^*,0}V^*$, where $T(V)\subset {\cal G}$ is a
left-symmetric algebra given by equation (6.21) and $(\rho^*,0)$ is
its bimodule since its sub-adjacent Lie algebra ${\cal G}(T(V))$ is
a Lie subalgebra of ${\cal G}$, and $T$ can be identified as an
element in $T(V)\otimes V^*\subset
(T(V)\ltimes_{\rho^*,0}V^*)\otimes (T(V)\ltimes_{\rho^*,0}V^*)$.

{\bf Proof}\quad Let $\{e_1,\cdots,e_n\}$ be a basis of ${\cal G}$. Let
$\{v_1,\cdots, v_m\}$ be a basis of $V$ and $\{ v_1^*,\cdots, v_m^*\}$ be
its dual basis. Set $T(v_i)=\sum\limits_{k=1}^na_{ik}e_k,
i=1,\cdots, m$. Then, we have
$$T=\sum_{i=1}^m T(v_i)\otimes v_i^*=\sum_{i=1}^m\sum_{k=1}^n
a_{ik}e_k\otimes v_i^*\in T(V)\otimes
V^*\subset (T(V)\ltimes_{\rho^*,0}V^*)\otimes (T(V)\ltimes_{\rho^*}V^*).$$
Therefore we have
\begin{eqnarray*}
-r_{12}\cdot r_{13} &=&-\sum_{i,j=1}^m\{T(v_i)\cdot T(v_j)\otimes v_i^*\otimes
v_j^*+v_i^*\otimes v_j^*\otimes T(\rho(T(v_j))v_i)\};\\
r_{12}\cdot r_{23} &=&\sum_{i,j=1}^m\{v_i^*\otimes T(v_i)\cdot T(v_j)\otimes
v_j^*-v_i^*\otimes v_j^*\otimes T(\rho(T(v_i))v_j)\};\\
\mbox{}[r_{13},r_{23}] &=&\sum_{i,j=1}^m \{
-v_i^*\otimes T(\rho(T(v_i))v_j)\otimes v_j^*
+T(\rho(T(v_j))v_i)\otimes v_j^*\otimes v_i^*\\
&\mbox{}&+ v_i^*\otimes v_j^*\otimes [T(v_i), T(v_j)]\}.
\end{eqnarray*}
Since $T$ is an ${\cal O}$-operator associated to $\rho$ and $T(u)\cdot T(v)=T(\rho (T(u)) v))$ for any $u,v\in V$, we know that
$r$ is a symmetric solution of $S$-equation in $T(V)\ltimes_{\rho^*,0}V^*$.\hfill $\Box$

{\bf Remark }\quad The above theorem also provides a method to
construct the symmetric solutions of $S$-equation (hence
parak\"ahler Lie algebras). We would like to point out that this
method starts from Lie algebras and only involves the representation
theory of Lie algebras (to find the ${\cal O}$-operators), which
avoids involving the nonassociativity of left-symmetric algebras.

{\bf Corollary 6.10}\quad Let $(A,\cdot)$ be a left-symmetric
algebra. Then
$$r=\sum_{i=1}^n (e_i\otimes e_i^*+e_i^*\otimes e_i)\eqno (6.28)$$
is a symmetric solution of $S$-equation in $A\ltimes_{L^*,0} A^*$,
where $\{e_1,\cdots, e_n\}$ is a basis of $A$ and $\{e_1^*,\cdots,
e_n^*\}$ is its dual basis.  Moreover, $r$ is nondegenerate and the
induced 2-cocycle ${\cal B}$ of $A\ltimes_{L^*,0} A^*$ is given by
$${\cal B}(x+a^*,y+b^*)=\langle r^{-1}(x+a^*),y+b^*\rangle  =\langle x,b^*\rangle  +\langle a^*,y\rangle  ,\;\;\forall x,y\in A, a^*,b^*\in A^*.\eqno (6.29)$$

{\bf Proof}\quad Let $V=A$, $\rho=L$ and $T=id$ in Theorem 6.9. Then
the conclusion follows immediately. \hfill $\Box$

{\bf Remark}\quad Comparing with Theorem 5.5, we know that (the
non-symmetric) $T=\sum\limits_{i=1}^ne_i\otimes e_i^*$ induces a
left-symmetric bialgebra structure on $A\ltimes_{{\rm ad}^*,-R^*}
A^*$, whereas the above (symmetric) $r=T+T^{21}$ induces a
left-symmetric bialgebra structure on $A\ltimes_{L^*,0} A^*$.

\section{Comparison between Lie bialgebras and left-symmetric
bialgebras}

In this section, we first recall some facts on Lie bialgebras from
Chapters 1-2 in  [CP] by Chari and Pressley. In fact, a Lie
bialgebra is the Lie algebra ${\cal G}$ of a Poisson-Lie group $G$
equipped with additional structures induced from the Poisson
structure on $G$ and a Poisson-Lie group is a Lie group with a
Poisson structure compatible with the group operation in a certain
sense. Poisson-Lie groups play an important role in symplectic
geometry and quantum group theory ([D],[KM],[LW]).

{\bf Definition 7.1}\quad Let ${\cal G}$ be a Lie algebra. A Lie
bialgebra structure on ${\cal G}$ is a skew-symmetric linear map
$\delta_{\cal G}: {\cal G}\rightarrow {\cal G}\otimes {\cal G}$ such
that $\delta_{\cal G}^*:{\cal G}^*\otimes {\cal G}^*\rightarrow
{\cal G}^*$ is a Lie bracket on ${\cal G}^*$ and $\delta$ is a
1-cocycle of ${\cal G}$ associated to ${\rm ad}\otimes 1+1\otimes
{\rm ad}$ with values in ${\cal G}\otimes {\cal G}$. We denote it by
$({\cal G},{\cal G}^*)$ or $({\cal G},\delta_{\cal G})$.

Let ${\cal G}$ be a Lie algebra. A bilinear form $B(\;,\;)$ on
${\cal G}$ is called invariant if
$$B([x,y],z)=B(x,[y,z]),\;\;\forall x,y,z\in {\cal G}.\eqno (7.1)$$

{\bf Definition 7.2}\quad A Manin triple is a triple of Lie algebras
$({\cal P}, {\cal P}_+,{\cal P}_-)$ together with a non-degenerate
symmetric invariant bilinear form $B(\;,\;)$ on ${\cal P}$ such that

(1) ${\cal P}_+$ and ${\cal P}_-$ are Lie subalgebras of ${\cal P}$;

(2) ${\cal P}={\cal P}_+\oplus {\cal P}_-$ as vector spaces;

(3) $B({\cal P}_+,{\cal P}_+)=B({\cal P}_-,{\cal P}_-)=0$.

\noindent Two Manin triple $({\cal P}_1,{\cal P}_{1,+},{\cal
P}_{1,-})$ and $({\cal P}_2,{\cal P}_{2,+},{\cal P}_{2,-})$ with the
bilinear forms $B_1(\;,\;)$ and $B(\;,\;)$ respectively are
isomorphic if there exists an isomorphism of Lie algebras
$\varphi:{\cal P}_1\rightarrow {\cal P}_2$ such that
$$\varphi({\cal P}_{1,+})={\cal P}_{2,+},\;\;\varphi({\cal P}_{1,-})={\cal P}_{2,-},\;
B_1(x,y)=B_2(\varphi(x),\varphi(y)),\;\;\forall x,y\in {\cal
P}_1.\eqno (7.2)$$

In particular, if there is a Lie algebra structure ${\cal G}\bowtie
{\cal G}^*$ on ${\cal G}\oplus {\cal G}^*$ such that ${\cal G}$ and
${\cal G}^*$ are Lie subalgebras and the natural symmetric bilinear
form (also see equation (6.29))
$$(x+a^*,y+b^*)=\langle a^*,y\rangle  +\langle x,b^*\rangle  ,\;\;\forall x,y\in {\cal G},a^*,b^*\in {\cal G}^*,\eqno (7.3)$$
is invariant, then $({\cal G}\bowtie {\cal G}^*, {\cal G}, {\cal
G}^*)$ is a (standard) Manin triple. It is known that every Manin
triple is isomorphic to such a standard Manin triple.

{\bf Theorem 7.1}\quad Let $({\cal G},[\;,\;]_{\cal G})$ be a Lie
algebra and $({\cal G}^*, [\;,\;]_{{\cal G}^*})$ be a Lie algebra
structure on its dual space ${\cal G}^*$. Then the following
conditions are equivalent:

(1) $({\cal G}\bowtie {\cal G}^*, {\cal G},{\cal G}^*)$ is a
standard Manin triple with the bilinear form (7.3);

(2) $({\cal G}, {\cal G}^*, {\rm ad}^*_{\cal G}, {\rm ad}^*_{{\cal
G}^*})$ is a matched pair of Lie algebras;

(3) $({\cal G},{\cal G}^*)$ is a Lie bialgebra.

{\bf Proposition 7.2}\quad Let $({\cal G}, {\cal G}^*)$ be a Lie
bialgebra. Then there is a canonical Lie bialgebra structure on
${\cal G}\oplus {\cal G}^*$ such that the inclusions $i_1:{\cal
G}\rightarrow {\cal G}\oplus {\cal G}^*$ and $i_2:{\cal
G}^*\rightarrow {\cal G}\oplus {\cal G}^*$ into the two summands are
homomorphisms of Lie bialgebras. Such a structure is called a
classical (Drinfeld) double of ${\cal G}$.

{\bf Definition 7.3}\quad A Lie bialgebra $({\cal G}, \delta)$ is
called coboundary if there exists a $r\in {\cal G}\otimes {\cal G}$
such that
$$\delta(x)=[x\otimes 1+1\otimes x, r],\;\;\forall x\in {\cal
G}.\eqno (7.4)$$

{\bf Theorem 7.3}\quad Let ${\cal G}$ be a Lie algebra and $r\in
{\cal G}\otimes {\cal G}$. Then the map $\delta: {\cal G}\rightarrow
{\cal G}\otimes {\cal G}$ defined by equation (7.4) induces a Lie
bialgebra structure on ${\cal G}$ if and only if the following two
conditions are satisfied:

(a) $[x\otimes 1+1\otimes x,\;r_{12}+r_{21}]=0$ for any $x\in {\cal
G}$;

(b) $[x\otimes 1\otimes 1+1\otimes x\otimes 1+1\otimes 1\otimes
x,\;\; [r_{12},r_{13}]+[r_{12},r_{23}]+[r_{13},r_{23}]]=0$ for any
$x\in {\cal G}$,

\noindent where the notations $r_{12}, r_{21}, r_{13}, r_{23}$ are
given as the notations after Proposition 5.3 for a Lie algebra.

{\bf Corollary 7.4}\quad Let ${\cal G}$ be a Lie algebra and $r\in
{\cal G}\otimes {\cal G}$. If $r$ is skew-symmetric and $r$
satisfies
$$[r_{12},r_{13}]+[r_{12},r_{23}]+[r_{13},r_{23}]=0,\eqno (7.5)$$
then the map $\delta: {\cal G}\rightarrow {\cal G}\otimes {\cal G}$
defined by equation (7.4) induces a Lie bialgebra structure on
${\cal G}$.

{\bf Definition 7.4}\quad Let ${\cal G}$ be a Lie algebra and $r\in
{\cal G}\otimes {\cal G}$. Equation (7.5) is called classical
Yang-Baxter equation (CYBE) in ${\cal G}$.

The facts above are from Chapters 1-2 in  [CP]. Next we recall some
results from the literature.

{\bf Proposition 7.5}\quad Let ${\cal G}$ be a Lie algebra and $r\in
{\cal G}\otimes {\cal G}$.

(1) [D]\quad Suppose $r$ is skew-symmetric and nondegenerate. Then
$r$ is a solution of CYBE in ${\cal G}$ if and only if the
isomorphism ${\cal G}^*\rightarrow {\cal G}$ induced by $r$,
regarded as a bilinear form on ${\cal G}$ is a 2-cocycle of ${\cal
G}$. Therefore under this situation, ${\cal G}$ is a symplectic Lie
algebra.

(2) ([Ku3])\quad Suppose $r$ is skew-symmetric. Then $r$ is a
solution of CYBE in ${\cal G}$ if and only if $r$ is an ${\cal
O}$-operator associated to ${\rm ad}^*$, that is, $r$ satisfies
$$[r(a^*),r(b^*)]=r({\rm ad}^*(r(a^*))b^*-{\rm
ad}^*(r(b^*))a^*),\;\;\forall a^*,b^*\in {\cal G}^*.\eqno (7.6)$$

{\bf Remark}\quad Let ${\cal G}$ be a Lie algebra equipped with a
nondegenerate symmetric invariant bilinear form. Let $r\in {\cal
G}\otimes {\cal G}$. Then $r$ can be identified as a linear map from
${\cal G}$ to ${\cal G}$. If $r$ is skew-symmetric, then $r$ is a
solution of CYBE in ${\cal G}$ if and only if $r$ satisfies equation
(2.7), that is, the operator form of CYBE.

{\bf Proposition 7.6} ([Bai3]) \quad Let ${\cal G}$ be a Lie algebra
and $\rho:{\cal G}\rightarrow gl(V)$ be a representation of ${\cal
G}$. Let $\rho^*:{\cal G}\rightarrow gl(V^*)$ be the dual
representation of $\rho$ and $T:V\rightarrow {\cal G}$ be a linear
map. Then
$$r=T-T^{21}\eqno (7.7)$$
is a skew-symmetric solution of CYBE in ${\cal
G}\ltimes_{\rho^*}V^*$ if and only if $T$ is an ${\cal O}$-operator
associated to $\rho$.

Therefore, roughly speaking, a symmetric solution of the
$S$-equation corresponds to the symmetric part of an ${\cal
O}$-operator, whereas a skew-symmetric solution of the classical
Yang-Baxter equation corresponds to the skew-symmetric part of an
${\cal O}$-operator.

{\bf Proposition 7.7} ([Bai3]) \quad Let $(A,\cdot)$ be a
left-symmetric algebra. Then
$$r=\sum_{i=1}^n (e_i\otimes e_i^*-e_i^*\otimes e_i)\eqno (7.7)$$
is a skew-symmetric solution of  CYBE in ${\cal G}(A)\ltimes_{L^*}
{\cal G}(A)^*$, where $\{e_1,\cdots, e_n\}$ is a basis of $A$ and
$\{e_1^*,\cdots, e_n^*\}$ is its dual basis.  Moreover, $r$ is
nondegenerate and the induced 2-cocycle ${\cal B}$ of ${\cal
G}(A)\ltimes_{L^*} {\cal G}(A)^*$ given by
$${\cal B}(x+a^*,y+b^*)=\langle r^{-1}(x+a^*),y+b^*\rangle  ,\;\;\forall x,y\in A, a^*,b^*\in A^*,\eqno (7.8)$$
is precisely equation (2.16).

The facts above and the results in the previous sections allow us to
compare left-symmetric bialgebras and Lie bialgebras  in terms of
the following properties:  structures on the corresponding Lie
groups, 1-cocycles of Lie algebras, matched pairs of Lie algebras,
Lie algebra structures on the direct sum of the Lie algebras in the
matched pairs, bilinear forms on the direct sum of the Lie algebras
in the matched pairs, double structures on the direct sum of the Lie
algebras in the matched pairs, algebraic equations associated to
coboundary cases, nondegenerate solutions and related geometric
interpretation, ${\cal O}$-operators and constructions from
left-symmetric algebras. We list them in Table 1. From this table,
we observe that there is a clear analogy between them and in
particular, parak\"ahler Lie groups correspond to Poisson-Lie groups
whose Lie algebras are Lie bialgebras in this sense.

Moreover, since classical Yang-Baxter equations can be regarded as
``classical limits" of quantum Yang-Baxter equations ([Be]), we
believe that there should exist an analogue (``quantum
$S$-equations" ) of the quantum Yang-Baxter equations. We use a
question mark in Table 1 to denote these still-to-be-found ``quantum
$S$-equations".

\begin{table}[t]\caption{Comparison between Lie bialgebras and left-symmetric bialgebras}
\begin{tabular}{|c|c|c|}
\hline Algebras & Lie bialgebras & Left-symmetric bialgebras\\\hline
Corresponding Lie groups & Poisson-Lie groups & parak\"ahler Lie
groups
\\\hline
1-cocycles of Lie algebras & $1\otimes {\rm ad}+{\rm ad}\otimes 1$ &
$L\otimes 1+1\otimes {\rm ad}$\\\hline Matched pairs of Lie algebras
& $({\cal G},{\cal G}^*, {\rm ad}^*_{\cal G}, {\rm ad}^*_{{\cal
G}^*})$ & $({\cal G}(A), {\cal G}(A^*), L^*_{A}, L^*_{A^*})$\\\hline
Lie algebra  structures on   & Manin triples & phase spaces\\
the direct sum of the Lie  &&\\ algebras in the matched
pairs&&\\\hline Bilinear forms on  & symmetric & skew-symmetric
\\\cline{2-3} the direct sum of the Lie & $\langle x+a^*,y+b^*\rangle  $ & $\langle x+a^*,y+b^*\rangle  $\\algebras in the
matched pairs &$=\langle x,b^*\rangle  +\langle a^*,y\rangle  $ &
$=-\langle x,b^*\rangle  +\langle a^*,y\rangle  $\\\cline{2-3}
&invariant & 2-cocycles\\\hline Double structures on  & Drinfeld
Doubles & symplectic doubles\\ the direct sum of the Lie&&\\algebras
in the matched pairs&&\\\hline Algebraic equations associated &
skew-symmetric solutions & symmetric solutions
\\\cline{2-3}
to coboundary cases& CYBEs in Lie algebras &$S$-equations in
left-symmetric\\&&algebras\\\hline Nondegenerate
solutions & 2-cocycles of Lie algebras & 2-cocycles of left-symmetric\\
&&algebras\\\cline{2-3} & symplectic structures & Hessian
structures\\\hline ${\cal O}$-operators & associated to ${\rm ad}^*$
& associated to $L^*$
\\\cline{2-3} &skew-symmetric parts & symmetric parts\\\hline
Constructions from & $r=\sum\limits_{i=1}^n (e_i\otimes
e_i^*-e_i^*\otimes e_i)$ & $r=\sum\limits_{i=1}^n (e_i\otimes
e_i^*+e_i^*\otimes e_i)$\\\cline{2-3} left-symmetric algebras &
induced bilinear
forms & induced bilinear forms\\
& $\langle x+a^*,y+b^*\rangle  $ & $\langle x+a^*,y+b^*\rangle  $\\
&$=-\langle x,b^*\rangle  +\langle a^*,y\rangle  $ & $=\langle
x,b^*\rangle  +\langle a^*,y\rangle  $\\\hline Quantum equations &
Quantum Yang-Baxter & ?\\& equations &\\\hline
\end{tabular}
\end{table}

At the end of this section, we consider the case that a
left-symmetric bialgebra is also a Lie bialgebra.

{\bf Theorem 7.7}\quad Let $(A,A^*,\alpha,\beta)$ be a
left-symmetric bialgebra. Then $({\cal G}(A),{\cal G}(A^*))$ is a
Lie bialgebra if and only if
$$\langle R_\cdot^*(x)a^*,R_\circ^*(b^*)y\rangle  +\langle R_\cdot^*(x)b^*,R_\circ^*(a^*)y\rangle
=\langle R_\cdot^*(y)b^*,R_\circ^*(a^*)x\rangle  +\langle
R_\cdot^*(y)a^*,R_\circ^*(b^*)x\rangle  , \eqno (7.9)$$ for any
$x,y\in A^*,a^*,b^*\in A^*$.

{\bf Proof}\quad Denote the left-symmetric product on $A$ by
$``\cdot''$ and the left-symmetric product on $A^*$ by $``\circ''$.
Then $({\cal G}(A),{\cal G}(A^*))$ is a Lie bialgebra if and only if
for any $x,y\in A$, $a^*,b^*\in A^*$,
$$\leqno(*)\;\;\; {\rm ad}_\cdot^*(x)[a^*,b^*]-[{\rm ad}_\cdot^*(x)a^*,b^*]
-[a^*,{\rm ad}_\cdot^*(x)b^*]+{\rm ad}_\cdot^*({\rm
ad}_\circ^*(a^*)x)b^*-{\rm ad}_\cdot^*({\rm
ad}_\circ^*(b^*)x)a^*=0;$$
$$\leqno(**)\;\;\; {\rm ad}_\circ^*(a^*)[x,y]-[{\rm ad}_\circ^*(a^*)x,y]
-[x,{\rm ad}_\circ^*(a^*)y]+{\rm ad}_\circ^*({\rm
ad}_\cdot^*(x)a^*)y-{\rm ad}_\circ^*({\rm ad}_\cdot^*(y)a^*)x=0.$$
Since ${\rm ad}_{\cdot}^*=L_\cdot^*-R_\cdot^*$ and $({\cal
G}(A),{\cal G}(A^*),L_\cdot^*,L_\circ^*)$ is a matched pair of Lie
algebras, the equation $(*)$ is reduced to
\begin{eqnarray*}
&&-R_\cdot^*(x)[a^*,b^*]+[R_\cdot^*(x)a^*,b^*]+[a^*,R_\cdot^*(x)b^*]
-R_\cdot^*({\rm ad}_{\circ}^*(a^*)x)b^*-L_\cdot^*(R_\circ^*(a^*)x)b^*\\
&&+R_\cdot^*({\rm
ad}_{\circ}^*(b^*)x)a^*+L_\cdot^*(R_\circ^*(b^*)x)a^*=0.
\end{eqnarray*}
Let the equation above act on $y\in A$ and note that $\langle
R_\cdot^*(x)a^*,y\rangle  =\langle L_\cdot^*(y)a^*,x\rangle  $ for
any $x,y\in A$ and $a^*\in A^*$, we obtain
\begin{eqnarray*}
&&-L_\cdot^*(y)[a^*,b^*]+L_\cdot^*({\rm ad}^*(b^*)y)a^*- L_\cdot^*({\rm ad}^*(a^*)y)b^*+[a,L_\cdot^*(y)b^*]+R_\cdot^*(y)b^*\circ a^*\\
&&-[b^*,L_\cdot^*(y)a^*]-R_\cdot^*(y)a^*\circ b^*=0.
\end{eqnarray*}
Using the condition that $({\cal G}(A),{\cal
G}(A^*),L_\cdot^*,L_\circ^*)$ is a matched pair of Lie algebras
again, we know that
$$-L_\cdot^*(R_\circ^*(b^*)y)a^*+ L_\cdot^*(R_\circ^*(a^*)y)b^*+R_\cdot^*(y)b^*\circ a^*-R_\cdot^*(y)a^*\circ b^*=0,$$
which gives equation (7.9) by acting on $x\in A$. Similarly, from
equation $(**)$, we can get
$$-L_\circ^*(R_\cdot^*(y)b^*)x+L_\circ^*(R_\cdot^*(x)b^*)y+R_\circ^*(b^*)y\cdot x-R_\circ^*(b^*)x\cdot y=0,$$
for any $x,y\in A$ and $b^*\in A^*$, which gives the same equation
(7.9) by acting on $a^*\in A^*$.\hfill $\Box$

{\bf Corollary 7.8}\quad Let $(A,\cdot)$ be a left-symmetric and
$r\in A\otimes A$ be a symmetric solution of $S$-equation in $A$.
Suppose the left-symmetric algebra structure on $A^*$ is induced by
$r$ from equation (6.4). Then there exists a Lie bialgebra structure
$({\cal G}(A),{\cal G}(A^*))$ if and only if
$$-r(L_\cdot^*(x)R_\cdot^*(y)a^*)+[x,r(R_\cdot^*(y)a^*]
-(y\cdot r(a^*)-r({\rm ad}_\cdot^*(y)a^*))\cdot x$$$$
=-r(L_\cdot^*(y)R_\cdot^*(x)a^*)+[y,r(R_\cdot^*(x)a^*] -(x\cdot
r(a^*)-r({\rm ad}_\cdot^*(x)a^*))\cdot y,\eqno (7.10)$$ for any
$x,y\in A$ and $a^*\in A^*$.

{\bf Corollary 7.9}\quad Let $(A,A^*,\alpha,\beta)$ be a
left-symmetric bialgebra. If equation (7.9) is satisfied, then there
are two Lie algebra structures ${\cal
G}(A)\bowtie_{L^*_\circ}^{L^*_\cdot}{\cal G}(A^*)$ and ${\cal
G}(A)\bowtie_{{\rm ad}^*_\circ}^{{\rm ad}^*_\cdot}{\cal G}(A^*)$ on
the direct sum $A\oplus A^*$ of the underlying vector spaces of $A$
and $A^*$ such that both ${\cal G}(A)$ and ${\cal G}(A^*)$ are Lie
subalgebras and the bilinear form given by equation (2.16) is a
2-cocycle of ${\cal G}(A)\bowtie_{L^*_\circ}^{L^*_\cdot}{\cal
G}(A^*)$ and the bilinear form (7.3) is invariant on ${\cal
G}(A)\bowtie_{{\rm ad}^*_\circ}^{{\rm ad}^*_\cdot}{\cal G}(A^*)$.

{\bf Remark}\quad Obviously the two Lie algebras above are not
isomorphic in general. On the other hand, if a Lie bialgebra $({\cal
G},{\cal G}^*)$ whose Lie algebra structure on ${\cal G}^*$ is
induced by a non-degenerate classical $r$-matrix ([D], [Se]), then
both ${\cal G}$ and ${\cal G}^*$ are symplectic Lie algebras.
Therefore there is a compatible left-symmetric algebra structure on
${\cal G}$ and ${\cal G}^*$ respectively.

{\bf Example 7.1}\quad Let $({\cal G},\omega)$ be a symplectic Lie
algebra. Then there is a Lie bialgebra whose Lie algebra structure
in ${\cal G}^*$ is given by a non-degenerate classical $r$-matrix as
follows (cf. [DiM]).
$$\delta(x)= [x\otimes 1+1\otimes x, r],\;\;\forall x\in {\cal G},\eqno (7.11)$$
where $r:{\cal G}^*\rightarrow {\cal G}$ is given by
$\omega(x,y)=\langle r^{-1}(x),y\rangle  $. On the other hand, there
exists a left-symmetric algebra structure $``\cdot"$ on ${\cal G}$
given by equation (2.14), that is, $\omega(x\cdot
y,z)=-\omega(y,[x,z])$ for any $x,y,z\in {\cal G}$. Moreover, there
exists a compatible left-symmetric algebra structure on the Lie
algebra ${\cal G}^*$ given by
$$a^*\circ b^*=r^{-1}(r(a^*)\cdot r(b^*)),\forall a^*,b^*\in {\cal G}^*.\eqno (7.12)$$
Furthermore, it is easy to know that
$$L_\cdot^*(x)a^*=r^{-1}[x,r(a^*)],\;R_\cdot^*(x)a^*=-r^{-1}(r(a^*)\cdot x)),\;
L_\circ^*(a^*)x=[r(a^*),x],\;R_\circ^*(a^*)x=-x\cdot r(a^*),\eqno
(7.13)$$ for any $x\in A,a^*\in A^*$. Therefore according to Theorem
2.5, $({\cal G},{\cal G}^*)$ (as left-symmetric algebras) is a
left-symmetric bialgebra if and only if ${\cal G}$ is 2-step
nilpotent, that is, $[[x,y],z]=0$ for any $x,y,z\in {\cal G}$. In
this case, we can know that it is equivalent to $[x,y]\cdot z=0$ for
any $x,y,z\in {\cal G}$. Therefore, equation (7.9) holds naturally.

\section*{ Acknowledgments}

The author thanks Professors P. Etingof, I.M. Gel'fand, B.A.
Kupershmidt, and C. Woodward for important suggestion and great
encouragement. He is very grateful of referees' important
suggestion. He also thanks Professors J. Lepowsky, Y.-Z. Huang and
H.S. Li for the hospitality extended to him during his stay at
Rutgers, The State University of New Jersey and for valuable
discussions. This work was supported in part by S.S. Chern
Foundation for Mathematical Research, the National Natural Science
Foundation of China (10571091, 10621101), NKBRPC (2006CB805905),
Program for New Century Excellent Talents in University and K.C.
Wong Education Foundation.

\end{document}